\documentclass[12pt]{article}

\usepackage{latexsym,amssymb,amsmath,ulem}
\usepackage{color}
\usepackage{graphicx}
\usepackage{latexsym,amssymb,amsmath,ulem}
\usepackage{amsfonts}
\usepackage{latexsym}
\usepackage{graphicx}
\usepackage{color}

\def\tto{\;{\lower 1pt \hbox{$\rightarrow$}}\kern -10pt
\hbox{\raise 2pt \hbox{$\rightarrow$}}\;}

\newtheorem{theorem}{Theorem}[section]

\newtheorem{corollary}{Corollary}[section]
\newtheorem{lemma}{Lemma}[section]
\newtheorem{remark}{Remark}[section]
\newtheorem{example}{Example}[section]
\newtheorem{definition}{Definition}[section]

\numberwithin{equation}{section}

\newcommand{\qed}{\hfill $\square$}
\newcommand{\QATOP}[2]{\genfrac{}{}{0pt}{}{#1}{#2}}
\normalsize

\begin{document}

\title{\textbf{Characterizations of robust and stable duality for linearly
perturbed uncertain optimization problems }}
\author{N. Dinh\thanks{
International University, Vietnam National University - HCMC, Linh Trung
ward, Thu Duc district, Ho Chi Minh city, Vietnam (ndinh@hcmiu.edu.vn). }, \
\ M.A. Goberna\thanks{%
Department of Mathematics, University of Alicante, Spain (mgoberna@ua.es){}}%
, \ \ M.A. Lopez\thanks{%
Department of Mathematics, University of Alicante, Spain, and CIAO,
Federation University, Ballarat, Australia (marco.antonio@ua.es) }, \ \ M.
Volle\thanks{%
Avignon University, LMA EA 2151, Avignon, France
(michel.volle@univ-avignon.fr)} }
\maketitle
\date{}


\textbf{Abstract.} We introduce a robust optimization model consisting in a
family of perturbation functions giving rise to certain pairs of dual
optimization problems in which the dual variable depends on the uncertainty
parameter. The interest of our approach is illustrated by some examples,
including uncertain conic optimization and infinite optimization via
discretization. The main results characterize desirable robust duality
relations (as robust zero-duality gap) by formulas involving the
epsilon-minima or the epsilon-subdifferentials of the objective function.
The two extreme cases, namely, the usual perturbational duality (without
uncertainty), and the duality for the supremum of functions (duality
parameter vanishing) are analyzed in detail.\medskip 

\medskip 


\section{Introduction}

In this paper we consider a \textit{family of perturbation functions} 
\begin{equation*}
F_{u}:X\times Y_{u}\rightarrow \mathbb{R}_{\infty }:=\mathbb{R}\cup
\{+\infty \},\text{ with }u\in U,
\end{equation*}%
and where $X$ and $Y_{u},$ $u\in U,$ are given locally convex Hausdorff
topological vector spaces (briefly, lcHtvs), the index set $U$ is called the 
\textit{uncertainty set }of the family\textit{, }$X$ is its \textit{decision
space}, and each $Y_{u}$ is a \textit{parameter space}. Note that our model
includes a parameter space $Y_{u},$ depending on $u\in U,$ which is a
novelty with respect to the \textquotedblleft classical" robust duality
scheme (see \cite{LJL11} and references therein, where a unique parameter
space $Y$ is considered), allowing us to cover a wider range of applications
including uncertain optimization problems under linear perturbations of the
objective function. The significance of our approach is illustrated along
the paper by relevant cases extracted from deterministic optimization with
linear perturbations, uncertain optimization without linear perturbations,
uncertain conic optimization and infinite optimization.\textbf{\ }The
antecedents of the paper are described in the paragraphs devoted to the
first two cases in Section 2.

We associate with each family $\left\{ F_{u}:u\in U\right\} $\ of
perturbation functions corresponding optimization problems whose definitions
involve continuous linear functionals on the decision and the parameter
spaces. We denote by $0_{X},$ ${0}_{_{X}}^{\ast },$ $0_{u},$ and $%
0_{u}^{\ast },$\textit{\ }the null vectors of $X,$ its topological dual $%
X^{\ast },$ $Y_{u},$ and its topological dual $Y_{u}^{\ast },$ respectively.
The optimal value of a minimization (maximization, respectively) problem $%
\mathrm{(}$P$\mathrm{)}$ is denoted by $\inf \mathrm{(}$P$\mathrm{)}$ ($\sup 
\mathrm{(}$P$\mathrm{)}$); in particular, we write $\min \mathrm{(}$P$%
\mathrm{)}$ ($\max \mathrm{(}$P$\mathrm{)}$) whenever the optimal value of $%
\mathrm{(}$P$\mathrm{)}$ is attained. We adopt the usual convention that $%
\inf \mathrm{(P)=+\infty }$ ($\sup \mathrm{(}$P$\mathrm{)=-\infty }$)\textbf{%
\ }when the problem\textbf{\ }$\mathrm{(P)}$ has no feasible solution.%
\textbf{\ }The associated optimization problems are the following:

\begin{itemize}
\item \textit{Linearly perturbed uncertain problems}: for each $(u,x^{\ast
})\in U\times X^{\ast },$%
\begin{equation*}
\mathrm{(}\text{\textrm{P}}_{u}\mathrm{)}{_{{x^{\ast }}}:}\quad \quad
\inf_{x\in X}\left\{ {{F_{u}}(x,{0_{u}})-\left\langle {x^{\ast },x}%
\right\rangle }\right\} .
\end{equation*}

\item \textit{Robust counterpart of }$\left\{ \mathrm{(}P_{u}\mathrm{)}{_{{%
x^{\ast }}}}\right\} _{u\in U}:$ 
\begin{equation*}
{\mathrm{(RP)}_{{x^{\ast }}}:}\quad \quad \inf_{x\in X}\left\{ {\mathop
{\sup }\limits_{u\in U}{F_{u}}(x,{0_{u}})-\left\langle {x^{\ast },x}%
\right\rangle }\right\} .
\end{equation*}
\end{itemize}

Denoting by $F_{u}^{\ast }:X^{\ast }\times Y_{u}^{\ast }\rightarrow 
\overline{\mathbb{R}}$, where $\overline{\mathbb{R}}:=\mathbb{R\cup \{\pm
\infty \}}$, the \textit{Fenchel conjugate} of $F_{u}$, namely, 
\begin{equation*}
F_{u}^{\ast }(x^{\ast },y_{u}^{\ast }):=\sup\limits_{(x,y_{u})\in X\times
Y_{u}}\Big\{\langle x^{\ast },x\rangle +\langle y_{u}^{\ast },y_{u}\rangle
-F_{u}(x,y_{u})\Big\},\quad (x^{\ast },y_{u}^{\ast })\in X^{\ast }\times
Y_{u}^{\ast },
\end{equation*}%
we now introduce the corresponding dual problems:

\begin{itemize}
\item \textit{Perturbational dual of }${(\text{\textrm{P}}_{u})_{{x^{\ast }}}%
}:$ 
\begin{equation*}
\mathrm{(}\text{\textrm{D}}_{u}\mathrm{)}{_{{x^{\ast }}}:}\quad \quad %
\mathop {\sup }\limits_{y_{u}^{\ast }\in Y_{u}^{\ast }}-F_{u}^{\ast }({%
x^{\ast }},y_{u}^{\ast }).
\end{equation*}%
Obviously,%
\begin{equation*}
\sup \mathrm{(}\text{\textrm{D}}_{u}\mathrm{)}_{x^{\ast }}\leq \inf \mathrm{(%
}\text{\textrm{P}}_{u}\mathrm{)}_{x^{\ast }}\leq \inf \mathrm{(}\text{%
\textrm{RP}}\mathrm{)}_{x^{\ast }},\forall u\in U.
\end{equation*}

\item \textit{Optimistic dual of} $\mathrm{(}$\textrm{RP}$\mathrm{)}%
_{x^{\ast }}:$ 
\begin{equation*}
\mathrm{(}\text{\textrm{ODP}}\mathrm{)}_{x^{\ast }}:\quad \quad
\sup_{(u,y_{u}^{\ast })\in \Delta }-F_{u}^{\ast }(x^{\ast },y_{u}^{\ast }),
\end{equation*}%
where 
$\Delta :=\left\{ {\left( {u,y_{u}^{\ast }}\right) :u\in U,\ {y^{\ast }}\in
Y_{u}^{\ast }}\right\} $ is the disjoint union of the spaces ${Y_{u}^{\ast }}
$. 
We have 
\begin{equation*}
\sup {(\mathrm{ODP})_{{x^{\ast }}}}=\mathop {\sup }\limits_{u\in U}{\left( {{%
\mathrm{D}_{u}}}\right) _{{x^{\ast }}}}\leq \inf {\mathrm{(RP)}_{{x^{\ast }}%
}.}
\end{equation*}
\end{itemize}

We are interested in the following desirable robust duality properties:

$\bullet $ \textit{Robust duality} is said to hold at $x^{\ast }$ if $\inf {%
\mathrm{(RP)}_{{x^{\ast }}}}=\sup {(\mathrm{ODP})_{{x^{\ast }}}}$,

$\bullet $ \textit{Strong robust duality} at $x^{\ast }$ means $\inf {%
\mathrm{(RP)}_{{x^{\ast }}}}=\max {(\mathrm{ODP})_{{x^{\ast }}}}$,

$\bullet $ \textit{Reverse strong robust duality }at $x^{\ast }$ means $\min 
{\mathrm{(RP)}_{{x^{\ast }}}}=\sup {(\mathrm{ODP})_{{x^{\ast }}}}$,

$\bullet $ \textit{Min-max robust duality} at $x^{\ast }$ means $\min {%
\mathrm{(RP)}_{{x^{\ast }}}}=\max {(\mathrm{ODP})_{{x^{\ast }}}}$.

Each of the above desirable properties is said to be \textit{stable} when it
holds for any $x^{\ast }\in X^{\ast }$. The main results of this paper
characterize these properties in terms of formulas involving the $%
\varepsilon $-minimizers and $\varepsilon $-subdifferentials of the
objective function of the robust counterpart problem \textrm{(RP)}$%
_{0_{X}^{\ast }}$, namely, the function%
\begin{equation*}
p:=\sup\limits_{u\in U}F_{u}(\cdot ,0_{u}).
\end{equation*}

Theorem \ref{thm31} characterizes robust duality at a given point $x^{\ast
}\in X^{\ast }$ as a formula for the inverse mapping of the $\varepsilon $%
-subdifferential at $x^{\ast }$ without any convexity assumption. The same
is done in Theorem \ref{thm41} to characterize strong robust duality. In the
case when a primal optimal solution does exist we give a formula for the
exact minimizers of $p-x^{\ast }$ to characterize dual strong (resp.
min-max) robust duality at $x^{\ast }$, see Theorem \ref{thm51} (resp.
Theorem \ref{thm52}). We show that stable robust duality gives rise\textbf{\ 
}to a formula for the $\varepsilon $-subdifferential of $p$ (Theorem \ref%
{thm61}, see also Theorem \ref{thm31}). The same is done for stable strong
robust duality (Theorem \ref{lem71}). A formula for the exact
subdifferential of $p$ is provided in relation with robust duality at
appropriate points (Theorem \ref{thm81}). The most simple possible formula
for the exact subdifferential of $p$ (the so-called \textit{Basic Robust
Qualification condition}) is studied in detail in Theorem \ref{thm82}. All
the results from Sections 1-8 are specified for the two extreme cases (the
case with no uncertainty and the one in absence of linear perturbations),
namely, Cases 1 and 2 in Section 2 (for the sake of brevity, we do not give
the specifications for Cases 3 and 4). It is worth noticing the generality
of the mentioned results (as they do not require any assumption on the
involved functions) and the absolute self containment of their proofs. The
use of convexity in the data will be addressed in a forthcoming paper.

\section{Special Cases and Applications}

In this section we make explicit the meaning of the robust duality of the
general model introduced in Section 1, composed by a family of perturbation
functions together with its corresponding optimization problems. We are
doing this by exploring the extreme case with no uncertainty, the extreme
case in absence of perturbations, and two other significant situations. In
all these cases, we propose \textit{ad hoc} families of perturbation
functions allowing to apply the duality results to given optimization
problems, either turning back to variants of well-known formulas for
conjugate functions or proposing\textbf{\ }new ones.

\noindent Let us recall the robust duality formula, $\inf {\mathrm{(RP)}_{{%
x^{\ast }}}=}\sup \mathrm{(}\mathrm{ODP}\mathrm{)}_{x^{\ast }},$ i.e., 
\begin{equation}
\;\;\mathop {\inf }\limits_{x\in X}\mathop {\sup }\limits_{u\in U}\left\{ {{%
F_{u}}\left( {x,{0_{u}}}\right) -\left\langle {{x^{\ast }},x}\right\rangle }%
\right\} =\mathop {\sup }\limits_{\left( {u,y_{^{u}}^{\ast }}\right) \in
\Delta }-{F_{u}^{\ast }}\left( {{x^{\ast }},{y_{u}^{\ast }}}\right) .
\label{eq1}
\end{equation}%
We firstly study the two extreme cases: the case with no uncertainty and the
one with no perturbations.

\textbf{Case 1. The case with no uncertainty:} Deterministic optimization
with linear perturbations deals with parametric problems of the form: 
\begin{equation*}
\mathrm{(P)}_{{x^{\ast }}}:\quad \quad {\mathop {\inf }\limits_{x\in X}}%
\left\{ {{f}(x)-\left\langle {{x^{\ast }},x}\right\rangle }\right\} {,}
\end{equation*}%
where $f:X\rightarrow \mathbb{R}_{\infty }$ (i.e., $f\in (\mathbb{R}_{\infty
})^{X}$) is the \textit{nominal objective function} and the \textit{parameter%
} is ${x^{\ast }\in X}^{\ast }.$ Taking a singleton uncertainty set $%
U=\left\{ u_{0}\right\} ,$\ $Y_{u_{0}}=Y$ and $F_{u_{0}}=F\ $such that ${%
F\left( {x,{0_{Y}}}\right) ={f}(x)}$ for all $x\in X,$ \eqref{eq1} reads 
\begin{equation}
\mathop {\inf }\limits_{x\in X}\left\{ {F\left( {x,{0_{Y}}}\right)
-\left\langle {{x^{\ast }},x}\right\rangle }\right\} =\mathop {\sup }%
\limits_{{y^{\ast }}\in {Y^{\ast }}}-{F^{\ast }}\left( {{x^{\ast }},{y^{\ast
}}}\right) ,  \label{11}
\end{equation}%
which is \textit{the fundamental perturbational duality formula} \cite%
{Bot-book1}, \cite{Rock74}, \cite{Za02}. Stable and strong robust duality
theorems are given in \cite{BuJeWu06} (see also \cite{DGLV17} and \cite{JL09}
for infinite optimization problems).

\textbf{Case 2. The case with no linear perturbations:} Uncertain
optimization without linear perturbations deals with families of problems of
the form

\begin{equation*}
\mathrm{(P):}\quad \quad \left\{ {\mathop {\inf }\limits_{x\in X}{f_{u}}(x)}%
\right\} _{u\in U},
\end{equation*}%
where $f_{u}\in (\mathbb{R}_{\infty })^{X},$ $u\in U.$ Taking $F\ $such that 
$F_{u}(x,0_{u})=f_{u}(x)$ for all $u\in U,$ the problem $\mathrm{(}$\textrm{P%
}$_{u}\mathrm{)}{_{{0}_{_{X}}^{\ast }}}$ represents here the scenario of $%
\mathrm{(P)}$ corresponding to $u\in U,$\ while ${\mathrm{(RP)}_{{0}%
_{_{X}}^{\ast }}}$ is the \textit{robust} \textit{counterpart} (also called 
\textit{pessimistic} or \textit{minmax counterpart} in the robust
optimization literature) of $\mathrm{(P),}$ namely, 
\begin{equation*}
\text{\textrm{(RP)}}\quad \quad \inf_{x\in X}\mathop {\sup }\limits_{u\in U}{%
{f_{u}}(x)}.
\end{equation*}%
\newline
For instance, if $F_{u}(x,y_{u})=f_{u}(x),\ $for all $y_{u}\in Y_{u},$ and $%
\operatorname{dom}f_{u}\neq \emptyset $, we have 
\begin{equation}
F_{u}^{\ast }\left( {x,y_{u}^{\ast }}\right) =\left\{ 
\begin{array}{ll}
{f_{u}^{\ast }\left( {x^{\ast }}\right) ,} & \mathrm{if}\;\;{y_{u}^{\ast }}%
=0_{u}^{\ast }, \\ 
{+\infty ,} & \mathrm{if}\;\;y_{u}^{\ast }\neq 0_{u}^{\ast }.%
\end{array}%
\right.   \label{12}
\end{equation}%
Then \eqref{eq1} writes 
\begin{equation}
{\left( {\mathop {\sup }\limits_{u\in U}{f_{u}}}\right) ^{\ast }}({x^{\ast }}%
)=\mathop {\inf }\limits_{u\in U}f_{u}^{\ast }({x^{\ast }}),  \label{13}
\end{equation}%
which amounts, for $x^{\ast }={0}_{_{X}}^{\ast },$ to the $\inf -\sup $%
\textit{\ duality in robust optimization,}\ also called \textit{robust
infimum} (recall that any constrained optimization problem can be reduced to
an unconstrained one by summing up the indicator function of the feasible
set to the objective function): 
\begin{equation*}
\inf_{x\in X}\sup_{u\in U}f_{u}(x)=\sup_{u\in U}\inf_{x\in X}f_{u}(x).
\end{equation*}%
Robust duality theorems without linear perturbations are given in \cite%
{WSS15} for a special class of uncertain non-convex optimization problems
while \cite{DGLV17} provides robust strong duality theorems for uncertain
convex optimization problems which are expressed in terms of the closedness
of suitable sets regarding the vertical axis of $X^{\ast }\times \mathbb{R}.$

\textbf{Case 3. Conic optimization problem with uncertain constraints:} 
\textit{Consider the uncertain problem} 
\begin{equation*}
\mathrm{(P):}\ \ \ \ \ \ \ \ \left\{ {\inf_{x\in X}f(x)\;}\;\text{\textrm{%
s.t.}}{\text{ }\;{H_{u}}(x)\in -{S_{u}}}\right\} _{u\in U},
\end{equation*}%
where, for each $u\in U$, $S_{u}$ is an ordering convex cone in ${Y}_{u},$ $%
H_{u}\colon X\rightarrow {Y}_{u}$, and $f\in (\mathbb{R}_{\infty })^{X}$.
Denote by ${S_{u}^{+}:=}\left\{ {y_{u}^{\ast }\in Y}_{u}^{\ast
}:\left\langle {y_{u}^{\ast },y_{u}}\right\rangle \geq 0,\forall {y_{u}\in }%
S_{u}\right\} $ the dual cone of $S_{u}$\textbf{.}

Problems of this type arise, for instance, in the production planning of
firms producing $n$\ commodities from uncertain amounts of resources by
means of technologies which depend on the available resources (e.g., the
technology differs when the energy is supplied by either fuel gas or a
liquid fuel). The problem associated with each parameter $u\in U$\ consists
of maximizing the cash-flow $c\left( x_{1},...,x_{n}\right) $\ of the total
production, with $x_{i}$\ denoting the production level of the $i$-th
commodity, $i=1,..,n.$\ The decision vector $x=\left( x_{1},...,x_{n}\right) 
$\ must satisfy a linear inequality system $A_{u}x\leq b_{u},$\ where the
matrix of technical coefficients $A_{u}$\ is $m_{u}\times n$\ and $b_{u}\in 
\mathbb{R}^{m_{u}},$\ for some $m_{u}\in \mathbb{N}.\ $Denoting by $\mathrm{i%
}_{\mathbb{R}_{+}^{n}}$ the indicator function of $\mathbb{R}_{+}^{n}$
(i.e., $\mathrm{i}_{\mathbb{R}_{+}^{n}}(x)=0,$ when $x\in \mathbb{R}%
_{+}^{n}, $ and $\mathrm{i}_{\mathbb{R}_{+}^{n}}(x)=+\infty ,$ otherwise),
the uncertain production planning problem can be formulated as%
\begin{equation*}
\mathrm{(P):}\ \ \ \ \ \ \ \ \left\{ {\inf_{x\in \mathbb{R}^{n}}f(x)=-c(x)+}%
\mathrm{i}_{\mathbb{R}_{+}^{n}}{(x)\;}\;\text{\textrm{s.t.}}{\text{ }%
\;A_{u}x-b_{u}\in -}\mathbb{R}_{+}^{m_{u}}\right\} _{u\in U},
\end{equation*}%
with the space $Y_{u}=\mathbb{R}^{m_{u}}$ depending on the uncertain
parameter $u.$

For each $u\in U$, define the perturbation function 
\begin{equation*}
{F_{u}}(x,{y_{u}})=\left\{ {%
\begin{array}{ll}
f{(x),} & {\mathrm{if}\;{H_{u}}(x)+{y_{u}}\in -{S_{u}}}, \\ 
{+\infty ,} & {\mathrm{else}}.%
\end{array}%
}\right.
\end{equation*}%
On the one hand, $\mathrm{(RP)}_{0_{X}^{\ast }}$ collapses to the robust
counterpart of $\mathrm{(P)}$\ in the sense of robust conic optimization
with uncertain constraints: 
\begin{equation*}
\mathrm{(RP):}\;\;\;\;\;\inf_{x\in X}f(x)\;\;\text{\textrm{s.t.}}\ \ \;{H_{u}%
}(x)\in -{S_{u}},\;\;\forall u\in U.
\end{equation*}%
On the other hand, it is easy to check that 
\begin{equation*}
F_{u}^{\ast }({x^{\ast }},{y}_{u}^{\ast })=\left\{ {%
\begin{array}{ll}
{{{\left( {f+y_{u}^{\ast }\circ {H_{u}}}\right) }^{\ast }}({x^{\ast }}),} & {%
\mathrm{if}\;y_{u}^{\ast }\in S_{u}^{+},} \\ 
{\ +\infty ,} & {\mathrm{else,}}%
\end{array}%
}\right.
\end{equation*}%
$\mathrm{(ODP)}_{0_{X}^{\ast }}$ is nothing else than the optimistic dual in
the sense of uncertain conic optimization:

\begin{equation*}
\mathrm{(ODP):}\quad \quad \sup\limits_{u\in U,{y}_{u}^{\ast }\in S_{u}^{+}}%
\mathop {\inf }\limits_{x\in X}\left\{ {f(x)+\left\langle {y_{u}^{\ast },{%
H_{u}}(x)}\right\rangle }\right\}
\end{equation*}%
(a special case when $Y_{u}=Y$, $S_{u}=S$ for all $u\in U$ is studied in 
\cite[page 1097]{DMVV} and \cite{LJL11}). Thus,

$\bullet$ \textit{Robust duality holds at $0_{X}^{\ast }$} means that $\inf 
\mathrm{(RP)}=\sup \mathrm{(}\mathrm{ODP}), $\ 

$\bullet $ \textit{Strong robust duality holds at $0_{X}^{\ast }$} means
that 
\begin{equation*}
\inf \left\{ {f(x):{H_{u}}(x)\in -{S_{u}},\forall u\in U}\right\}
=\max\limits_{\QATOP{u\in U\hfill }{{y_{u}^{\ast }}\in S_{u}^{+}\hfill }%
}\inf\limits_{x\in X}\left\{ {f(x)+\left\langle {y_{u}^{\ast },{H_{u}}(x)}%
\right\rangle }\right\} .
\end{equation*}%
Conditions for having such an equality are provided in \cite[Theorem 6.3]%
{DMVV}, \cite[Corollaries 5, 6]{DL17}, for the particular case $Y_{u}=Y$ for
all $u\in U$.

\textit{Strong robust duality and uncertain Farkas lemma:} We focus again on
the case where $Y_{u}=Y$ and $S_{u}=S$ for all $u\in U$. For a given $r\in 
\mathbb{R}$, let us\textbf{\ }consider the following statements:

\begin{itemize}
\item[(i)] $H_{u}(x)\in -S,\;\forall u\in U\quad \Longrightarrow \quad
f(x)\geq r$,

\item[(ii)] $\exists u\in U,\exists {y_{u}^{\ast }}\in S^{+}$ such that $%
f(x)+\left\langle {y_{u}^{\ast }},H_{u}(x)\right\rangle \geq r,\;\forall
x\in X.$
\end{itemize}

\noindent Then, it is true that the Strong robust duality holds at $%
0_{X}^{\ast }$ if and only if $\left[ (i)\Longleftrightarrow (ii)\right] $
for each $r\in \mathbb{R},$ which can be seen as an uncertain Farkas lemma.
For details\textbf{\ }see \cite[Theorem 3.2]{DMVV} (also \cite[Corollary 5
and Theorem 1]{DL17}).\textbf{\ }

It is worth noticing that when return to problem $\mathrm{(P)}$, a given
robust feasible solution $\overline{x}$\ is a minimizer if and only if $f(%
\overline{x})\leq f(x)$ for any robust feasible solution $x.$ So, a robust
(uncertain) Farkas lemma (with $r=f(\bar{x})$) will lead automatically to an
optimality test for $\mathrm{(P).}$ Robust conic optimization problems are
studied in \cite{BS06} and \cite{Vera17}.

\textbf{\textbf{Case 4. Discretizing infinite optimization problems:}} Let $%
f\in (\mathbb{R}_{\infty })^{X}$ and $g_{t}\in \mathbb{R}^{X}\;$for all$%
\;t\in T$ (a possibly infinite index set). Consider the set $U$ of nonempty
finite subsets of $T,$ interpreted as admissible perturbations of $T,$ and
the parametric optimization problem 
\begin{equation*}
\mathrm{(P):}\quad \left\{ \inf_{x\in X}f(x)\;\mathrm{s.t.}\text{ }%
g_{t}(x)\leq 0,\text{ }\forall t\in S\right\} _{S\in U}.
\end{equation*}%
Consider the parameter space $Y_{s}:={\mathbb{R}^{S}}$ (depending on $S$)
and the perturbation function ${F_{S}}:X\times {\mathbb{R}^{S}}\rightarrow 
\mathbb{R}_{\infty }$ such that, for any $x\in X$ and ${{{{\mu :=}({\mu _{s}}%
)}_{s\in S}\in }\mathbb{R}^{S},}$ 
\begin{equation*}
{F_{S}}\left( {x,{\mu }}\right) =\left\{ {%
\begin{array}{ll}
{f(x),} & {\mathrm{if}\;{g_{s}}(x)\leq -{\mu _{s}},\;\forall s\in S,} \\ 
{+\infty ,} & {\mathrm{else}}.%
\end{array}%
}\right.
\end{equation*}%
We now interpret the problems associated with the family of function
perturbations $\left\{ {F_{S}:S\in U}\right\} .$ One has $Y_{s}^{\ast }={%
\mathbb{R}^{S}}$ and 
\begin{equation*}
F_{S}^{\ast }\left( {{x^{\ast }},{\lambda }}\right) =\left\{ {%
\begin{array}{ll}
{{\left( {f+\sum\limits_{s\in S}{{\lambda _{s}g_{s}}}}\right) ^{\ast }}({%
x^{\ast }}),} & {\mathrm{if}\;{\lambda }\in \mathbb{R}_{+}^{S},} \\ 
{+\infty ,} & {\mathrm{else}}.%
\end{array}%
}\right.
\end{equation*}%
The robust counterpart at $0_{X}^{\ast },$ 
\begin{equation*}
\mathrm{(RP)}_{0_{X}^{\ast }}:\quad \inf f(x)\;\;\mathrm{s.t.}\ \;\;{g_{t}}%
(x)\leq 0\text{\ for all }t\in T,
\end{equation*}%
is a general infinite optimization problem while the optimistic dual at $%
0_{X}^{\ast }$ is 
\begin{equation*}
\mathrm{(ODP)}_{0_{X}^{\ast }}:\quad \sup\limits_{S\in U,{{\lambda }\in 
\mathbb{R}_{+}^{S}}}\left\{ \mathop {\inf }\limits_{x\in X}\left( {%
f(x)+\sum\limits_{s\in S}{{\lambda _{s}g_{s}}(x)}}\right) \right\} ,
\end{equation*}%
or, equivalently, the Lagrange dual of $\mathrm{(RP)}_{0_{X}^{\ast }},$
i.e., 
\begin{equation*}
\mathrm{(ODP)}_{0_{{X}^{\ast }}}:\quad \quad \sup\limits_{{\lambda }\in 
\mathbb{R}_{+}^{(T)}}\left\{ \mathop {\inf }\limits_{x\in X}\left( {%
f(x)+\sum\limits_{t\in T}{{\lambda _{t}g_{t}}(x)}}\right) \right\} ,
\end{equation*}%
where, for each $\lambda =\left( {\lambda _{t}}\right) _{t\in T}\in \mathbb{R%
}_{+}^{(T)}$ (the subspace of $\mathbb{R}^{T}$ formed by the functions $%
\lambda $ whose support, supp$\lambda :=\left\{ t\in T:{{\lambda _{t}\neq 0}}%
\right\} ,$ is finite), 
\begin{equation*}
\sum\limits_{t\in T}{{\lambda _{t}g_{t}}(x):}=\left\{ 
\begin{array}{ll}
\sum\limits_{t\in \text{supp}\lambda }{{\lambda _{t}g_{t}}(x),} & {\mathrm{if%
}\;\lambda \neq 0,} \\ 
0, & \mathrm{if}\;\lambda =0.%
\end{array}%
\right.
\end{equation*}%
Following \cite[Section 8.3]{GL98}, we say that $\mathrm{(RP)}_{0_{X}^{\ast
}}$ is \textit{discretizable} if there exists a sequence $\left(
S_{r}\right) _{r\in \mathbb{N}}\subset U$ such that 
\begin{equation}
\inf \mathrm{(RP)}_{0_{X}^{\ast }}=\lim_{r}\inf \left\{ f(x):{g_{t}}(x)\leq
0,\ \forall t\in S_{r}\right\} ,  \label{2.10}
\end{equation}%
and it is \textit{reducible }if there exists $S\in U$ such that 
\begin{equation*}
\inf \mathrm{(RP)}_{0_{X}^{\ast }}=\inf \left\{ f(x):{g_{t}}(x)\leq 0,\
\forall t\in S\right\} .
\end{equation*}%
Obviously, $\inf \mathrm{(RP)}_{0_{X}^{\ast }}=-\infty $ entails that $%
\mathrm{(RP)}_{0_{X}^{\ast }}$\ is reducible which, in turn, implies that $%
\mathrm{(RP)}_{0_{X}^{\ast }}$ is discretizable.

Discretizable and reducible problems are important in practice. Indeed, on
the one hand, discretization methods generate sequences $\left( S_{r}\right)
_{r\in \mathbb{N}}\subset U$ satisfying \ref{2.10} when $\mathrm{(RP)}%
_{0_{X}^{\ast }}$ is discretizable; discretization methods for linear and
nonlinear semi-infinite programs have been reviewed in \cite[Subsection 2.3]%
{GL17} and \cite{LS07}, while a hard infinite optimization problem has been
recently solved via discretization in \cite{LR17}. On the other hand,
replacing the robust counterpart (a hard semi-infinite program when the
uncertainty set is infinite) of a given uncertainty optimization problem,
when it is reducible, by a finite subproblem allows many times to get the
desired tractable reformulation (see e.g., \cite{BEN09} and \cite{BJL13}).

\begin{example}[Discretizing linear infinite optimization problems]
\label{Example1}Consider the problems introduced in Case 4 above, with $%
f(\cdot ):=\left\langle c^{\ast },\cdot \right\rangle $ and ${g_{t}}%
(x):=\left\langle a_{t}^{\ast },\cdot \right\rangle -b_{t},$ where $c^{\ast
},a_{t}^{\ast }\in X^{\ast }$ and $b_{t}\in \mathbb{R},$ for all $t\in T.$
Then, $\mathrm{(RP)}_{0_{X}^{\ast }}$ collapses to the linear infinite
programming problem%
\begin{equation*}
\mathrm{(RP)}_{0_{X}^{\ast }}:\;\;\;\inf \quad \left\langle c^{\ast
},x\right\rangle \;\;\mathrm{s.t.}\ \ \left\langle a_{t}^{\ast
},x\right\rangle \leq b_{t},\ \forall t\in T,
\end{equation*}%
whose feasible set we denote by $A.$ So, $\inf \mathrm{(RP)}_{0_{X}^{\ast
}}=\inf_{x\in X}\left\{ \left\langle c^{\ast },x\right\rangle +i_{A}\left(
x\right) \right\} .\;$ We assume that $A\neq \emptyset .$ \newline
Given $S\in U$ and ${{{\mu ,\lambda }\in }}\mathbb{R}^{S},$ 
\begin{equation}
{F_{S}}\left( {x,{\mu }}\right) =\left\{ {%
\begin{array}{ll}
\left\langle c^{\ast },x\right\rangle {,} & {\mathrm{if}\;\ \left\langle
a_{s}^{\ast },x\right\rangle \leq b_{s}-{\mu _{s}},\;\forall s\in S,} \\ 
{+\infty ,} & {\mathrm{else,}}%
\end{array}%
}\right.  \label{2.11}
\end{equation}%
and 
\begin{equation}
F_{S}^{\ast }\left( {{x^{\ast }},{\lambda }}\right) =\left\{ {%
\begin{array}{ll}
{{{\sum\limits_{s\in S}{{\lambda _{s}b_{s}}}}},} & {\mathrm{if}\;{{{%
\sum\limits_{s\in S}{{\lambda _{s}a_{s}^{\ast }}}}}}}=x^{\ast }{{-c^{\ast }}}%
\text{ }\mathrm{and}\text{ }{{\lambda _{s}\geq 0\,,\ }\forall }s\in S, \\ 
{+\infty ,} & {\mathrm{else.}}%
\end{array}%
}\right.  \label{2.12}
\end{equation}%
Hence, $\mathrm{(ODP)}_{0_{X}^{\ast }}$ collapses to the so-called \textit{%
Haar dual problem \cite{GLV14}} of $\mathrm{(RP)}_{0_{X}^{\ast }},$ 
\begin{equation*}
\mathrm{(ODP)}_{0_{X}^{\ast }}:\quad \sup \left\{ -{{{\sum\limits_{t\in 
\operatorname{supp}{\lambda }}{{\lambda _{t}b_{t}:-{{{\sum\limits_{t\in \operatorname{%
supp}{\lambda }}}}\lambda _{t}a_{t}^{\ast }=}}}}c^{\ast },}}\text{ }{\lambda 
}\in \mathbb{R}_{+}^{(T)}\right\} ,
\end{equation*}%
i.e., 
\begin{equation}
\sup \mathrm{(ODP)}_{0_{X}^{\ast }}=-\inf_{S\in U,{{{\lambda }\in }}\mathbb{R%
}_{+}^{S}}\left\{ {{{\sum\limits_{s\in S}{{\lambda _{s}b_{s}:}%
\sum\limits_{s\in S}{{\lambda _{s}a_{s}^{\ast }=}}}}-c^{\ast }}}\right\} .
\label{2.13}
\end{equation}%
From (\ref{2.13}), if $\inf \mathrm{(RP)}_{0_{X}^{\ast }}=\max \mathrm{(ODP)}%
_{0_{X}^{\ast }}\in \mathbb{R},$ then there exist $S\in U$ and ${{{\lambda }%
\in }}\mathbb{R}_{+}^{S}$ such that 
\begin{equation}
{{{\sum\limits_{s\in S}{{\lambda _{s}}}}}}\left( {{{{{{a_{s}^{\ast },}b_{s}}}%
}}}\right) {=-}\left( c^{\ast },\inf \mathrm{(RP)}_{0_{X}^{\ast }}\right) .
\label{2.14}
\end{equation}%
Let $A_{S}:=\left\{ x\in X:{\ \left\langle a_{s}^{\ast },x\right\rangle \leq
b_{s},\;\forall s\in S}\right\} .$ Given $x\in A_{S},$ from (\ref{2.14}), 
\begin{equation*}
0\geq {{{\sum\limits_{s\in S}{{\lambda _{s}}}}}}\left( {\left\langle
a_{s}^{\ast },x\right\rangle -{{{{b_{s}}}}}}\right) {=-}\left\langle c^{\ast
},x\right\rangle +\inf \mathrm{(RP)}_{0_{X}^{\ast }}.
\end{equation*}%
Since 
\begin{equation*}
\inf \mathrm{(RP)}_{0_{X}^{\ast }}\leq \left\langle c^{\ast },x\right\rangle
,\forall x\in A_{S},
\end{equation*}%
\begin{equation}
\inf \mathrm{(RP)}_{0_{X}^{\ast }}=\inf \left\{ \left\langle c^{\ast
},x\right\rangle :{\left\langle a_{s}^{\ast },x\right\rangle \leq
b_{s},\;\forall s\in S}\right\} ,  \label{2.15}
\end{equation}%
so that $\mathrm{(RP)}_{0_{X}^{\ast }}$ is reducible. Conversely, if (\ref%
{2.15}) holds with $\inf \mathrm{(RP)}_{0_{X}^{\ast }}\in \mathbb{R}$ and $%
\operatorname{cone}\left\{ \left( {{{{{{a_{t}^{\ast },}b_{t}}}}}}\right) :t\in
T\right\} +\mathbb{R}_{+}\left( 0_{X}^{\ast },1\right) $ is weak$^{\ast }$%
\textbf{-}closed, since $\inf \mathrm{(RP)}_{0_{X}^{\ast }}\leq \left\langle
c^{\ast },x\right\rangle $ is consequence of $\left\{ {\left\langle
a_{s}^{\ast },x\right\rangle \leq b_{s},\;\forall s\in S}\right\} ,$\ by the
non-homogeneous Farkas lemma in lcHtvs \cite{Chu66} and the closedness
assumption, there exist ${{{\lambda }\in }}\mathbb{R}_{+}^{S}$ and $\mu {\in 
}\mathbb{R}_{+}$ such that 
\begin{equation*}
{-}\left( c^{\ast },\inf \mathrm{(RP)}_{0_{X}^{\ast }}\right) ={{{%
\sum\limits_{s\in S}{{\lambda _{s}}}}}}\left( {{{{{{a_{s}^{\ast },}b_{s}}}}}}%
\right) {+}\mu \left( 0_{X}^{\ast },1\right) ,
\end{equation*}%
which implies that $\mu =0$ and $\inf \mathrm{(RP)}_{0_{X}^{\ast }}=\max 
\mathrm{(ODP)}_{0_{X}^{\ast }}{.}$ The closedness assumption holds when $X$
is finite dimensional (guaranteeing that any finitely generated convex cone
in $X^{\ast }\times \mathbb{R}$ is closed). So, as proved in \cite[Theorem
8.3]{GL98}, a linear semi-infinite program $\mathrm{(RP)}_{0_{X}^{\ast }}$
is reducible if and only if (\ref{2.15}) holds if and only if $\inf \mathrm{%
(RP)}_{0_{X}^{\ast }}=\max \mathrm{(ODP)}_{0_{X}^{\ast }}.$\newline
We now assume that $\inf \mathrm{(RP)}_{0_{X}^{\ast }}=\sup \mathrm{(ODP)}%
_{0_{X}^{\ast }}\in \mathbb{R}.$ By (\ref{2.13}), there exist sequences $%
\left( S_{r}\right) _{r\in \mathbb{N}}\subset U$ and $\ \left( {\lambda }%
_{r}\right) _{r\in \mathbb{N}},$ with ${\lambda }^{r}{\in }\mathbb{R}%
_{+}^{S_{r}}$ for all $r\in \mathbb{N},$\ such that%
\begin{equation*}
\lim_{r}\inf_{{\lambda }^{r}{\in }\mathbb{R}_{+}^{S_{r}}}\left\{ {{{%
\sum\limits_{s\in S_{r}}{{\lambda _{s}^{r}b_{s}:}\sum\limits_{s\in S_{r}}{{%
\lambda _{s}^{r}a_{s}^{\ast }=}}}}-c^{\ast }}}\right\} =-\sup \mathrm{(ODP)}%
_{0_{X}^{\ast }}.
\end{equation*}%
Denote $v_{r}:=-{{{\sum\limits_{s\in S_{r}}{{\lambda _{s}^{r}b_{s}.}}}}}$
Then, 
\begin{equation}
{{{\sum\limits_{s\in S_{r}}{{\lambda _{s}}}}}}\left( {{{{{{a_{s}^{\ast },}%
b_{s}}}}}}\right) {=-}\left( c^{\ast },v_{r}\right) ,  \label{2.16}
\end{equation}%
with $\lim_{r}v_{r}=\inf \mathrm{(RP)}_{0_{X}^{\ast }}.$ Let $A_{r}:=\left\{
x\in X:{\ \left\langle a_{s}^{\ast },x\right\rangle \leq b_{s},\;\forall
s\in S}_{r}\right\} ,$ $r\in \mathbb{N}.$ Given $x\in A_{r},$ from (\ref%
{2.16}), 
\begin{equation*}
0\geq {{{\sum\limits_{s\in S_{r}}{{\lambda _{s}^{r}}}}}}\left( {\left\langle
a_{s}^{\ast },x\right\rangle -{{{{b_{s}}}}}}\right) ={-}\left\langle c^{\ast
},x\right\rangle +v_{r}.
\end{equation*}%
Since $v_{r}\leq \left\langle c^{\ast },x\right\rangle $ for all $x\in
A_{r}, $ 
\begin{equation*}
v_{r}\leq \inf \left\{ \left\langle c^{\ast },x\right\rangle :{\left\langle
a_{s}^{\ast },x\right\rangle \leq b_{s},\;\forall s\in S}_{r}\right\} \leq
\inf \mathrm{(RP)}_{0_{X}^{\ast }}.
\end{equation*}%
Thus, 
\begin{equation*}
\lim_{r}\inf \left\{ \left\langle c^{\ast },x\right\rangle :{\left\langle
a_{s}^{\ast },x\right\rangle \leq b_{s},\;\forall s\in S}_{r}\right\} =\inf 
\mathrm{(RP)}_{0_{X}^{\ast }},
\end{equation*}%
i.e., $\mathrm{(RP)}_{0_{X}^{\ast }}$ is discretizable. Once again, the
converse is true in linear semi-infinite programming \cite[Corollary 8.2.1]%
{GL98}, but not in linear infinite programming.
\end{example}

\section{Robust Conjugate Duality}

We now turn back to the general perturbation function ${F_{u}}\colon X\times 
{Y_{u}}\rightarrow \mathbb{R}_{\infty },$ $u\in U$ and let $\Delta :=\left\{
(u,y_{u}^{\ast }):u\in U,y_{u}^{\ast }\in Y_{u}^{\ast }\right\} $ be the
disjoint union of the spaces $Y_{u}^{\ast }$. Recall that 
\begin{equation}
{\mathrm{(RP)}_{{x^{\ast }}}:}\quad \inf_{x\in X}\left\{ {\mathop {\sup }%
\limits_{u\in U}{F_{u}}\left( {x,{0_{u}}}\right) -\left\langle {{x^{\ast }},x%
}\right\rangle }\right\} ,  \label{RPx-star}
\end{equation}%
\begin{equation}
\mathrm{(ODP)}{_{{x^{\ast }}}:}\quad \ \ \sup\limits_{(u,\mathbf{y}%
_{u}^{\ast })\in \Delta }-F_{u}^{\ast }\left( {{x^{\ast }},y_{u}^{\ast }}%
\right) .  \label{ODPx-star}
\end{equation}

Define $p\in \overline{\mathbb{R}}^{X}$ and $q\in \overline{\mathbb{R}}%
^{X^\ast}$ such that 
\begin{equation}
p:=\mathop {\sup }\limits_{u\in U}{F_{u}}(\cdot ,{0_{u}})\ \text{\ \ \ and \
\ \ }q:=\mathop {\inf }\limits_{\left( {u,y_{u}^{\ast }}\right) \in \Delta
}\;F_{u}^{\ast }(\cdot ,y_{u}^{\ast }).  \label{pq}
\end{equation}%
One then has 
\begin{equation}
\left\{ {%
\begin{array}{l}
{{p^{\ast }}({x^{\ast }})=-\inf {{\mathrm{(RP)}}_{{x^{\ast }}}},\quad q({%
x^{\ast }})=-\sup \mathrm{(}{{\mathrm{ODP)}}_{{x^{\ast }}}\bigskip }} \\ 
{{q^{\ast }}=\mathop {\sup }\limits_{\left( {u,y_{u}^{\ast }}\right) \in
\Delta }{{\left( {F_{u}^{\ast }\left( \cdot {,y_{u}^{\ast }}\right) }\right) 
}^{\ast }}=\mathop {\sup }\limits_{u\in U}F_{u}^{\ast \ast }\left( \cdot {,{%
0_{u}}}\right) \leq p},%
\end{array}%
}\right.  \label{30}
\end{equation}%
and hence,

\begin{itemize}
\item \textit{Weak robust duality} always holds: 
\begin{equation}
p^{\ast }(x^{\ast })\leq q^{\ast \ast }(x^{\ast })\leq q(x^{\ast }),\text{
for all }x^{\ast }\in X^{\ast }.  \label{WD}
\end{equation}

\item \textit{Robust duality at }$x^{\ast }$ means: 
\begin{equation}
p^{\ast }(x^{\ast })=q(x^{\ast }).  \label{SD}
\end{equation}
\end{itemize}

Robust duality at $x^{\ast }$ also holds when either $p^{\ast }(x^{\ast
})=+\infty $ or $q(x^{\ast })=-\infty .$

As an illustration, consider Case 4 with linear data, as in Example \ref%
{Example1}. Then, $p\left( x\right) =\left\langle c^{\ast },x\right\rangle +%
\mathrm{i}_{A}\left( x\right) ,$ $\operatorname{dom}p=A,$ and so 
\begin{equation*}
p^{\ast }\left( 0_{X}^{\ast }\right) =\sup_{x\in \mathbb{R}^{n}}\left(
-p\left( x\right) \right) =-\inf_{x\in \mathbb{R}^{n}}\left\{ \left\langle
c^{\ast },x\right\rangle +\mathrm{i}_{A}\left( x\right) \right\} =-\inf 
\mathrm{(RP)}_{0_{X}^{\ast }}.
\end{equation*}%
Similarly, from (\ref{2.12}), 
\begin{equation*}
q\left( {{x^{\ast }}}\right) =\inf_{S\in U,{{{\lambda }\in }}\mathbb{R}%
^{S}}\left\{ {{{\sum\limits_{s\in S}{{\lambda _{s}b_{s}:}\sum\limits_{s\in S}%
{{\lambda _{s}a_{s}^{\ast }=}}}}x^{\ast }-c^{\ast }}}\right\} ,
\end{equation*}%
$\operatorname{dom}q=c^{\ast }+\operatorname{cone}\left\{ a_{t}^{\ast }:t\in T\right\} 
$ and 
\begin{equation}
q\left( 0_{X}^{\ast }\right) =\inf_{S\in U,{{{\lambda }\in }}\mathbb{R}%
_{+}^{S}}\left\{ {{{\sum\limits_{s\in S}{{\lambda _{s}b_{s}:}%
\sum\limits_{s\in S}{{\lambda _{s}a_{s}^{\ast }=}}}}-c^{\ast }}}\right\}
=-\sup \mathrm{(ODP)}_{0_{X}^{\ast }}.  \label{36}
\end{equation}

\subsection{Basic lemmas}

Let us introduce the necessary notations. Given a lcHtvs $Z$, an extended
real-valued function $h\in \overline{\mathbb{R}}^{Z}$, and $\varepsilon \in 
\mathbb{R}_{+}$, the set of $\varepsilon $-minimizers of $h$ is defined by 
\begin{equation*}
\varepsilon -\text{argmin }h:=\left\{ 
\begin{array}{ll}
\{z\in Z\,:\,h(z)\leq \inf_{Z}h+\varepsilon \}, & \mathrm{if}%
\;\;\inf\limits_{Z}h\in \mathbb{R}, \\ 
\emptyset , & \mathrm{if}\;\;\inf\limits_{Z}h\not\in \mathbb{R},%
\end{array}%
\right. 
\end{equation*}%
or, equivalently, 
\begin{equation*}
\varepsilon -\text{argmin }h=\{z\in h^{-1}(\mathbb{R})\,:\,h(z)\leq
\inf_{Z}h+\varepsilon \}.
\end{equation*}%
Note that $\varepsilon -$argmin $h\neq \emptyset $ when $\inf_{Z}h\in 
\mathbb{R}$ and $\varepsilon >0.$ Various calculus rules involving $%
\varepsilon -$argmin have been given in \cite{Volle95}.

The $\varepsilon -$subdifferential of $h$ at a point $a\in Z$ is the set
(see, for instance, \cite{HMSV95}) 
\begin{eqnarray*}
\partial ^{\varepsilon }h(a) &:=&\left\{ 
\begin{array}{ll}
\{z^{\ast }\in Z^{\ast }\,:\,h(z)\geq h(a)+\langle z^{\ast },z-a\rangle
-\varepsilon ,\forall z\in Z\}, & \mathrm{if}\;\;h(a)\in \mathbb{R}, \\ 
\emptyset , & \mathrm{if}\;\;h(a)\not\in \mathbb{R}%
\end{array}%
\right. \\
&=&\Big\{z^{\ast }\in (h^{\ast })^{-1}(\mathbb{R})\,:\,h^{\ast }(z^{\ast
})+h(a)\leq \langle z^{\ast },a\rangle +\varepsilon \Big\}.
\end{eqnarray*}

It can be checked that if $h\in \overline{\mathbb{R}}^{X}$ is convex and $%
h(a)\in \mathbb{R}$, then $\partial ^{\varepsilon }h(a)\neq \emptyset $ for
all $\varepsilon >0$ if and only if $h$ is lower semi-continuous at $a$.

The inverse of the set-valued mapping $\partial ^{\varepsilon
}h:Z\rightrightarrows Z^{\ast }$ is denoted by $M^{\varepsilon }h:Z^{\ast
}\rightrightarrows Z.$ For each $(\varepsilon ,z^{\ast })\in \mathbb{R}%
_{+}\times Z^{\ast }$, we have

\textbf{\ } 
\begin{equation*}
\Big(\partial ^{\varepsilon }h\Big)^{-1}(z^{\ast })=\Big(M^{\varepsilon }h%
\Big)(z^{\ast })=\varepsilon -\text{argmin }(h-z^{\ast }).
\end{equation*}%
Denoting by $\partial ^{\varepsilon }h^{\ast }(z^{\ast })$ the $\varepsilon $%
-subdifferential of $h^{\ast }$ at $z^{\ast }\in Z^{\ast }$, namely, 
\begin{equation*}
\partial ^{\varepsilon }h^{\ast }(z^{\ast })=\Big\{z\in (h^{\ast \ast
})^{-1}(\mathbb{R})\,:\,h^{\ast \ast }(z)+h^{\ast }(z^{\ast })\leq \langle
z^{\ast },z\rangle +\varepsilon \Big\},
\end{equation*}%
where $h^{\ast \ast }(z):=\sup\limits_{z^{\ast }\in Z^{\ast }}\{\langle
z^{\ast },z\rangle -h^{\ast }(z^{\ast })\}$ is the biconjugate of $h$, we
have 
\begin{equation*}
(M^{\varepsilon }h)(z^{\ast })\subset (\partial ^{\varepsilon }h^{\ast }\big)%
(z^{\ast }),\ \forall (\varepsilon ,z^{\ast })\in \mathbb{R}_{+}\times
Z^{\ast },
\end{equation*}%
with equality if and only if $h=h^{\ast \ast }.$ 

For each $\varepsilon \in \mathbb{R}_{+}$, we consider the set-valued
mapping $S^{\varepsilon }:X^{\ast }\rightrightarrows X$ as follows: 
\begin{equation*}
\begin{array}{ll}
S^{\varepsilon }(x^{\ast }) & :=\left\{ x\in p^{-1}(\mathbb{R}%
)\,:\,p(x)-\langle x^{\ast },x\rangle \leq -q(x^{\ast })+\varepsilon
\right\} .%
\end{array}
\label{J-ep}
\end{equation*}
If $q(x^{\ast })=-\infty $, then $S^{\varepsilon }(x^{\ast })=p^{-1}(\mathbb{%
R}).$ If $q(x^{\ast })=+\infty ,$ then $S^{\varepsilon }(x^{\ast
})=\emptyset .$

Since $p^{\ast }\leq q$, it is clear that 
\begin{equation}
S^{\varepsilon }(x^{\ast })\subset (M^{\varepsilon }p)(x^{\ast }),\ \forall
\varepsilon \geq 0,\ \forall x^{\ast }\in X^{\ast }.  \label{31}
\end{equation}

\begin{lemma}
\label{lem31} Assume that $\operatorname{dom} p\neq \emptyset $. Then, for each $%
x^{\ast }\in X^{\ast } $, the next statements are equivalent:

$\mathrm{(i)}$ \textit{Robust duality holds at }$x^{\ast }$\thinspace $,$
i.e., $p^{\ast }(x^{\ast })=q(x^{\ast })$,

$\mathrm{(ii)} $ $\left( M^{\varepsilon }p\right) (x^{\ast })=S^{\varepsilon
}(x^{\ast }),\quad \forall \varepsilon \geq 0$,

$\mathrm{(iii)}$ $\exists \bar{\varepsilon}>0:\left( M^{\varepsilon
}p\right) (x^{\ast })=S^{\varepsilon }(x^{\ast }),\quad \forall \varepsilon
\in ]0,\bar{\varepsilon}[$.
\end{lemma}

\noindent \textit{Proof.} $[\mathrm{(i)}\Rightarrow \mathrm{(ii)}]$ By
definition%
\begin{eqnarray*}
\left( M^{\varepsilon }p\right) (x^{\ast }) &=&\varepsilon -\operatorname{argmin}%
(p-x^{\ast }) \\
&=&\left\{ x\in p^{-1}(\mathbb{R})\,:\,p(x)-\langle x^{\ast },x\rangle \leq
-p^{\ast }(x^{\ast })+\varepsilon \right\} .
\end{eqnarray*}%
By $\mathrm{(i)}$ we thus have\ $\left( M^{\varepsilon }p\right) (x^{\ast
})=S^{\varepsilon }(x^{\ast }).$

$[\mathrm{(ii)}\Rightarrow \mathrm{(iii)}]$ It is obviously true.

$[\mathrm{(iii)}\Rightarrow \mathrm{(i)}]$ Since $p^{\ast }(x^{\ast })\leq
q(x^{\ast })$, $\mathrm{(i)}$ holds if $p^{\ast }(x^{\ast })=+\infty .$
Moreover, since $\operatorname{dom}p\neq \emptyset $, one has $p^{\ast }(x^{\ast
})\neq -\infty .$ Let now $p^{\ast }(x^{\ast })\in \mathbb{R}$. In order to
get a contradiction, assume that $p^{\ast }(x^{\ast })\not=q(x^{\ast })$.
Then $p^{\ast }(x^{\ast })<q(x^{\ast })$ and there exists $\varepsilon \in
]0,\bar{\varepsilon}[$ such that $p^{\ast }(x^{\ast })+\varepsilon
<q(x^{\ast })$. Since $\inf_{x\in X}\left\{ p(x)-\left\langle x^{\ast
},x\right\rangle \right\} =-p^{\ast }(x^{\ast })\in \mathbb{R}$ and $%
\varepsilon >0,$ we have $\varepsilon -\operatorname{argmin}(p-x^{\ast })\neq
\emptyset .$ Let us pick $x\in (M^{\varepsilon }p)(x^{\ast })=\varepsilon -%
\operatorname{argmin}(p-x^{\ast }).$ By $\mathrm{(iii)}$, we have $x\in
S^{\varepsilon }(x^{\ast })$ and 
\begin{equation*}
-p^{\ast }(x^{\ast })\leq p(x)-\langle x^{\ast },x\rangle \leq -q(x^{\ast
})+\varepsilon ,
\end{equation*}%
which contradicts $p^{\ast }(x^{\ast })+\varepsilon <q(x^{\ast })$.\qed

For each $\varepsilon \in \mathbb{R}_{+}$, let us introduce now the the
following set-valued mapping $J^{\varepsilon }:U\rightrightarrows X$: 
\begin{equation}  \label{J-ep}
J^{\varepsilon }(u):=\left\{ x\in p^{-1}(\mathbb{R})\,:\,p(x)\leq
F_{u}(x,0_{u})+\varepsilon \right\} ,
\end{equation}%
with the aim of making explicit the set $S^{\varepsilon }(x^{\ast })$. To
this purpose, given $\varepsilon _{1},\varepsilon _{2}\in \mathbb{R}_{+}$, $%
u\in U$, and $y_{u}^{\ast }\in Y_{u}^{\ast }$, let us introduce the
set-valued mapping $A_{(u,y_{u}^{\ast })}^{(\varepsilon _{1},\varepsilon
_{2})}:X^{\ast }\rightrightarrows X$ such that 
\begin{equation*}
A_{(u,y_{u}^{\ast })}^{(\varepsilon _{1},\varepsilon _{2})}(x^{\ast }):=%
\Big\{x\in J^{\varepsilon _{1}}(u)\,:\,(x,0_{u})\in (M^{\varepsilon
_{2}}F_{u})(x^{\ast },y_{u}^{\ast })\Big\}.
\end{equation*}

\begin{lemma}
\label{lem32} For each $x^{\ast }\in X^{\ast }$, $\varepsilon
_{1},\varepsilon _{2}\in \mathbb{R}_{+}$, $u\in U$, and $y_{u}^{\ast }\in
Y_{u}^{\ast }$, one has 
\begin{equation*}
A_{(u,y_{u}^{\ast })}^{(\varepsilon _{1},\varepsilon _{2})}(x^{\ast })\
\subset \ S^{\varepsilon _{1}+\varepsilon _{2}}(x^{\ast }).
\end{equation*}
\end{lemma}

\noindent \textit{Proof.} Let $x\in J^{\varepsilon _{1}}(u)$ be such that $%
(x,0_{u})\in (M^{\varepsilon _{2}}F_{u})(x^{\ast },y_{u}^{\ast }).$ Then we
have\textbf{\ }$F_{u}^{\ast }(x^{\ast },y_{u}^{\ast })\in \mathbb{R}$ and%
\textbf{\ }$F_{u}(x,0_{u})\in \mathbb{R}$. Moreover%
\begin{equation*}
F_{u}(x,0_{u})+\varepsilon _{1}\geq p(x)\geq F_{u}(x,0_{u})\in \mathbb{R}%
\text{,}
\end{equation*}%
implying $p(x)\in \mathbb{R}$ and, by (\ref{30}), 
\begin{eqnarray*}
p(x)-\langle x^{\ast },x\rangle &\leq &F_{u}(x,0_{u})-\langle x^{\ast
},x\rangle +\varepsilon _{1}\leq -F_{u}^{\ast }(x^{\ast },y_{u}^{\ast
})+\varepsilon _{1}+\varepsilon _{2} \\
&\leq &-q(x^{\ast })+\varepsilon _{1}+\varepsilon _{2},
\end{eqnarray*}%
that means $x\in S^{\varepsilon _{1}+\varepsilon _{2}}(x^{\ast })$ . \qed

\begin{lemma}
\label{lem33} Assume that 
\begin{equation}
\operatorname{dom}F_{u}\neq \emptyset ,\text{ }\forall u\in U.  \label{32}
\end{equation}%
Then, for each $x^{\ast }\in X^{\ast },\varepsilon \in \mathbb{R}_{+},\eta
>0 $, one has 
\begin{equation*}
{S^{\varepsilon }}({x^{\ast }})\ \ \subset \ \ \bigcup\limits_{\QATOP{u\in
U\hfill }{y_{u}^{\ast }\in Y_{u}^{\ast }\hfill }}\bigcup\limits_{\QATOP{%
\scriptstyle{\varepsilon _{1}}+{\varepsilon _{2}}=\varepsilon +\eta \hfill }{%
\scriptstyle{\varepsilon _{1}}\geq 0,\;{\varepsilon _{2}}\geq 0\hfill }%
}A_{(u,y_{u}^{\ast })}^{(\varepsilon _{1},\varepsilon _{2})}(x^{\ast }).
\end{equation*}
\end{lemma}

\noindent \textit{Proof.} Let $x\in p^{-1}(\mathbb{R})$ be such that $x\in
S^{\varepsilon }(x^{\ast })$, i.e., 
\begin{equation*}
p(x)-\langle x^{\ast },x\rangle \leq -q(x^{\ast })+\varepsilon .
\end{equation*}%
We then have, for any $\eta >0$, 
\begin{equation*}
q(x^{\ast })<\langle x^{\ast },x\rangle -p(x)+\varepsilon +\eta
\end{equation*}%
and, by definition of $q$ and $p,$ there exist $u\in U$, $y_{u}^{\ast }\in
Y_{u}^{\ast }$ such that 
\begin{equation}
F_{u}^{\ast }(x^{\ast },y_{u}^{\ast })\leq \langle x^{\ast },x\rangle
-p(x)+\varepsilon +\eta \leq \langle x^{\ast },x\rangle
-F_{u}(x,0_{u})+\varepsilon +\eta .  \label{33}
\end{equation}%
Since $p(x)\in \mathbb{R}$, $F_{u}^{\ast }(x^{\ast },y_{u}^{\ast })\neq
+\infty $. In fact, by \eqref{32}, $F_{u}^{\ast }(x^{\ast },y_{u}^{\ast
})\in \mathbb{R}$. Similarly, $F_{u}(x,0_{u})\in \mathbb{R}$. Setting 
\begin{equation*}
\alpha _{1}:=p(x)-F_{u}(x,0_{u}),\ \alpha _{2}:=F_{u}^{\ast }(x^{\ast
},y_{u}^{\ast })+F_{u}(x,0_{u})-\langle x^{\ast },x\rangle ,
\end{equation*}%
we get $\alpha _{1}\in \mathbb{R}_{+}$, $\alpha _{2}\in \mathbb{R}.$
Actually {$\alpha _{2}\geq 0$ since, by definition of conjugate, } 
\begin{equation*}
F_{u}^{\ast }(x^{\ast },y_{u}^{\ast })=\sup_{z\in X,y_{u}\in Y_{u}}\left\{
\langle x^{\ast },z\rangle +\langle y_{u}^{\ast },y_{u}\rangle
-F_{u}(z,y_{u})\right\} ,
\end{equation*}%
i.e., if $z=x$ and $y_{u}=0_{u},$%
\begin{equation*}
F_{u}^{\ast }(x^{\ast },y_{u}^{\ast })\geq \langle x^{\ast },x\rangle
-F_{u}(x,0_{u}),
\end{equation*}%
so that 
\begin{equation*}
F_{u}^{\ast }(x^{\ast },y_{u}^{\ast })+F_{u}(x,0_{u})-\langle x^{\ast
},x\rangle \geq 0.
\end{equation*}%
Then, by \eqref{33}, $0\leq \alpha _{1}+\alpha _{2}\leq \varepsilon +\eta $.
Consequently, there exist $\varepsilon _{1},\varepsilon _{2}\in \mathbb{R}%
_{+}$ such that $\alpha _{1}\leq \varepsilon _{1}$, $\alpha _{2}\leq
\varepsilon _{2}$, $\varepsilon _{1}+\varepsilon _{2}=\varepsilon +\eta $.
Now $\alpha _{1}\leq \varepsilon _{1}$ means that $x\in J^{\varepsilon
_{1}}(u)$ and $\alpha _{2}\leq \varepsilon _{2}$ means that $(x,0_{u})\in
(M^{\varepsilon _{2}}F_{u})(x^{\ast },y_{u}^{\ast })$, and we have $x\in
A_{(u,y_{u}^{\ast })}^{(\varepsilon _{1},\varepsilon _{2})}(x^{\ast })$. 
\qed

\medskip

For each $x^{\ast }\in X^{\ast }$, $\varepsilon \in \mathbb{R}_{+}$, let us
define 
\begin{eqnarray*}
\mathcal{A}^{\varepsilon }(x^{\ast }):= &&\bigcap\limits_{\eta
>0}\bigcup\limits_{\QATOP{u\in U\hfill }{y_{u}^{\ast }\in Y_{u}^{\ast
}\hfill }}\bigcup\limits_{\QATOP{\scriptstyle{\varepsilon _{1}}+{\varepsilon
_{2}}=\varepsilon +\eta \hfill }{\scriptstyle{\varepsilon _{1}}\geq 0,\;{%
\varepsilon _{2}}\geq 0\hfill }}A_{(u,y_{u}^{\ast })}^{(\varepsilon
_{1},\varepsilon _{2})}(x^{\ast }) \\
&=&\bigcap_{\eta >0}\bigcup\limits_{\QATOP{u\in U\hfill }{y_{u}^{\ast }\in
Y_{u}^{\ast }\hfill }}{{\bigcup\limits_{\QATOP{\scriptstyle{\varepsilon _{1}}%
+{\varepsilon _{2}}=\varepsilon +\eta \hfill }{\scriptstyle{\varepsilon _{1}}%
\geq 0,\;{\varepsilon _{2}}\geq 0\hfill }}}}\Big\{x\in J^{\varepsilon
_{1}}(u)\,:\,(x,0_{u})\in (M^{\varepsilon _{2}}F_{u})(x^{\ast },y_{u}^{\ast
})\Big\}
\end{eqnarray*}

\subsection{Robust duality}

We now can state the main result on characterizations of the robust
conjugate duality.

\begin{theorem}[Robust duality]
\label{thm31} Assume that $\operatorname{dom} p\neq \emptyset $. Then for each $%
x^{\ast }\in X^{\ast } $, the next statements are equivalent:

$\mathrm{(i)}$ $\inf {\mathrm{(RP)}_{{x^{\ast }}}}=\sup {(\mathrm{ODP})_{{%
x^{\ast }}}}$,

$\mathrm{( ii)} $ $\left( M^{\varepsilon }p\right) (x^{\ast })=\mathcal{A}%
^{\varepsilon }(x^{\ast }),\quad \forall \varepsilon \geq 0$,

$\mathrm{( iii)} $ $\exists \bar{\varepsilon}>0: \ \left( M^{\varepsilon
}p\right) (x^{\ast })=\mathcal{A}^{\varepsilon }(x^{\ast }),\quad \forall
\varepsilon \in ]0,\bar{\varepsilon}[$.
\end{theorem}

\noindent \textit{Proof.} We firstly claim that if $\operatorname{dom} p\neq
\emptyset $ then for each $x^* \in X^\ast$, $\varepsilon \in \mathbb{R}_+$,
it holds: 
\begin{equation}  \label{eqlem34}
S^\varepsilon (x^*) \ = \ \mathcal{A}^\varepsilon (x^*).
\end{equation}
In deed, as $\operatorname{dom} p\neq \emptyset $, \eqref{32} holds. It then
follows from Lemma \ref{lem33}, $S^\varepsilon (x^*) \ \subset \ \mathcal{A}%
^\varepsilon (x^*)$. On the other hand, for each $\eta > 0$, one has, by
Lemma \ref{lem32}, 
\begin{equation*}
\bigcup\limits_{\QATOP{u\in U\hfill }{y_{u}^{\ast }\in Y_{u}^{\ast }\hfill }%
} \bigcup\limits_{\QATOP{\scriptstyle{\varepsilon _{1}}+{\varepsilon _{2}}%
=\varepsilon +\eta \hfill }{\scriptstyle{\varepsilon _{1}}\geq 0,\;{%
\varepsilon _{2}}\geq 0\hfill }} A^{(\varepsilon_1, \varepsilon_2)}_{(u,
y_u^*)}(x^*) \ \subset \ S^{\varepsilon + \eta} (x^*).
\end{equation*}
Taking the intersection over all $\eta > 0$ we get 
\begin{equation*}
\mathcal{A}^\varepsilon (x^*) \subset \bigcap\limits_{\eta> 0}
S^{\varepsilon + \eta} (x^*) = S^\varepsilon (x^*),
\end{equation*}
and \eqref{eqlem34} follows. Taking into account the fact that (i) means $%
p^{\ast }(x^{\ast })=q(x^{\ast })$, the conclusions now follows from %
\eqref{eqlem34} and Lemma \ref{lem31}. \qed

\medskip

For the deterministic optimization problem with linear perturbations (i.e.,
non-uncertain case where $U$ is a singleton), the next result is a direct
consequence of Theorem \ref{thm31}.

\begin{corollary}[Robust duality for Case 1]
\label{corol31} Let $F\colon X\times Y\rightarrow \mathbb{R}_{\infty }$ be
such that $\operatorname{dom}F(\cdot ,0_{Y})\neq \emptyset $. Then, for each $%
x^{\ast }\in X^{\ast },$ the fundamental duality formula (\ref{11}) holds,
i.e., 
\begin{equation*}
\mathop {\inf }\limits_{x\in X}\left\{ {F(x,{0_{Y}})-\left\langle {{x^{\ast }%
},x}\right\rangle }\right\} =\mathop {\sup }\limits_{y\in {Y^{\ast }}%
}-F^{\ast }(x^{\ast },y^{\ast }),
\end{equation*}%
if and only any of the (equivalent) conditions (ii) or (iii) in Theorem \ref%
{thm31}\ holds, where 
\begin{equation}
\mathcal{A}^{\varepsilon }(x^{\ast })=\bigcap\limits_{\eta
>0}\bigcup\limits_{y^{\ast }\in Y^{\ast }}\Big\{x\in X\,:\,\left(
x,0_{Y}\right) \in \left( M^{\varepsilon +\eta }F\right) (x^{\ast },y^{\ast
})\Big\}.  \label{A-ep1}
\end{equation}
\end{corollary}

\noindent \textit{Proof.} Let $F_{u}=F:X\times Y\rightarrow \mathbb{R}%
_{\infty },\;p=F(\cdot ,0_{Y})$. In this case, one has, 
\begin{equation*}
J^{\varepsilon }(u)=\left\{ x\in X\,:\,F(x,0_{Y})\in \mathbb{R}\right\} ,\
\forall \varepsilon \geq 0,
\end{equation*}%
and $\mathcal{A}^{\varepsilon }(x^{\ast })$ will take the form \eqref{A-ep1}%
. The conclusion follows from Theorem \ref{thm31}. \qed

For uncertain optimization problem without linear perturbations, the
following result is a consequence of Theorem \ref{thm31}.

\begin{corollary}[Robust duality for Case 2]
\label{corol32} Let $(f_{u})_{u\in U}\subset \mathbb{R}_{\infty }^{X}$ be a
family of extended real-valued functions, $p=\sup_{u\in U}f_{u}$ be such
that $\operatorname{dom}p\neq \emptyset $. Then, for each $x^{\ast }\in X^{\ast
}, $ the $\inf -\sup $ duality in robust optimization (\ref{13}) holds,
i.e., 
\begin{equation*}
{\left( {\mathop {\sup }\limits_{u\in U}{f_{u}}}\right) ^{\ast }}({x^{\ast }}%
)=\mathop {\inf }\limits_{u\in U}f_{u}^{\ast }({x^{\ast }}),
\end{equation*}%
if and only any of the (equivalent) conditions (ii) or (iii) in Theorem \ref%
{thm31}\ holds, where 
\begin{equation}
{{\mathcal{A}^{\varepsilon }}({x^{\ast }})=\bigcap\limits_{\eta >0}{%
\;\bigcup\limits_{u\in U}{\bigcup\limits_{\QATOP{\scriptstyle{\varepsilon
_{1}}+{\varepsilon _{2}}=\varepsilon +\eta \hfill }{\scriptstyle{\varepsilon
_{1}}\geq 0,\;{\varepsilon _{2}}\geq 0\hfill }}}}}\left\{ {{{{{J^{{%
\varepsilon _{1}}}}(u)\cap ({M^{{\varepsilon _{2}}}}{f_{u}})({x^{\ast }})}}}}%
\right\} ,  \label{A-ep-xstar}
\end{equation}%
with%
\begin{equation*}
{{{{{J^{{\varepsilon _{1}}}}(u)=\{x\in p}}}}^{-1}(\mathbb{R}):\ f_{u}(x)\geq
p(x)-\varepsilon _{1}\}.
\end{equation*}
\end{corollary}

\noindent \textit{Proof.} Let $F_{u}(x,y_{u})=f_{u}(x),$ for all $u\in U$
and let $p=\mathop{\sup }\limits_{u\in U}{f_{u}}$. Then, by \eqref{J-ep},

\begin{equation}
{J^{\varepsilon }}(u)=\left\{ x\in p^{-1}(\mathbb{R})\,:\,f_{u}(x)\geq
p(x)-\varepsilon \right\} ,\ \forall \varepsilon \geq 0.  \label{34}
\end{equation}%
Moreover, recalling (\ref{12}),\ for each $u\in U$ such that $\operatorname{dom}%
f_{u}\neq \emptyset $, $(x^{\ast },y_{u}^{\ast })\in X^{\ast }\times
Y_{u}^{\ast }$, and $\varepsilon \geq 0$, 
\begin{equation}
{\left( {{M^{\varepsilon }}{F_{u}}}\right) \left( {{x^{\ast }},y_{u}^{\ast }}%
\right) =\left\{ {%
\begin{array}{ll}
{\left( {{M^{\varepsilon }}{f_{u}}}\right) \left( {x^{\ast }}\right) ,} & {%
\mathrm{if}\ \ y_{u}^{\ast }=0_{u}^{\ast },} \\ 
\emptyset , & \text{else}.%
\end{array}%
\;}\right. }  \label{35}
\end{equation}%
Finally, for each $(x^{\ast },\varepsilon )\in X^{\ast }\times \mathbb{R}%
_{+} $, $\mathcal{A}^{\varepsilon }({x^{\ast }})$ takes the form as in %
\eqref{A-ep-xstar}. The conclusion now follows from Theorem \ref{thm31}. 
\qed

\section{Strong Robust Duality}

We retain the notations in Section 3 and consider the robust problem $%
\mathrm{(RP)}_{x^{\ast }}$ and its robust dual problem $\mathrm{(ODP)}%
_{x^{\ast }}$ given in \eqref{RPx-star} and \eqref{ODPx-star}, respectively.
Let $p$ and $q$ be the functions defined by \eqref{pq} and recall the
relations in \eqref{30}, that is, 
\begin{equation*}
\left\{ {%
\begin{array}{l}
{{p^{\ast }}({x^{\ast }})=-\inf {{\mathrm{(RP)}}_{{x^{\ast }}}},\quad q({%
x^{\ast }})=-\sup \mathrm{(}{{\mathrm{ODP)}}_{{x^{\ast }}}\medskip }} \\ 
{{q^{\ast }}=\mathop {\sup }\limits_{\left( {u,y_{u}^{\ast }}\right) \in
\Delta }{{\left( {F_{u}^{\ast }\left( \cdot {,y_{u}^{\ast }}\right) }\right) 
}^{\ast }}=\mathop {\sup }\limits_{u\in U}F_{u}^{\ast \ast }\left( \cdot {,{%
0_{u}}}\right) \leq p}.%
\end{array}%
}\right.
\end{equation*}

In this section we will establish characterizations of \textit{strong robust
duality at $x^{\ast }$}. Recall that the strong robust duality holds at $x^*$
means that $\inf {\mathrm{(RP)}_{{x^{\ast }}}}=\max {(\mathrm{ODP})_{{%
x^{\ast }}}}$, which is the same as: 
\begin{equation*}
\exists (u,y_{u}^{\ast })\in \Delta :p^{\ast }(x^{\ast })=F_{u}^{\ast
}(x^{\ast },y_{u}^{\ast }).
\end{equation*}
For this, we need a technical lemma, but firstly, given $x^* \in X^\ast$, $\
u\in U$, $y_{u}^{\ast }\in Y_{u}^{\ast }$, and $\varepsilon \geq 0$, let us
introduce the set 
\begin{equation*}
B_{(u,y_{u}^{\ast })}^{\varepsilon }({x^{\ast }})= \bigcup\limits_{\QATOP{%
\scriptstyle{\varepsilon _{1}}+{\varepsilon _{2}} =\varepsilon \hfill }{%
\scriptstyle{\varepsilon _{1}}\geqslant 0,{\varepsilon _{2}}\geqslant
0\hfill }} A^{(\varepsilon_1, \varepsilon_2)}_{(u, y_u^*)}(x^*)
=\bigcup\limits_{\QATOP{\scriptstyle{\varepsilon _{1}}+{\varepsilon _{2}}%
=\varepsilon \hfill }{\scriptstyle{\varepsilon _{1}}\geqslant 0,{\varepsilon
_{2}}\geqslant 0\hfill }}\Big\{x\in J^{\varepsilon _1}(u)\, :\, (x,0_{u})
\in ( M^{ \varepsilon _{2}} F_{u}) (x^\ast ,y_u^\ast )\Big\}.
\end{equation*}

\begin{lemma}
\label{lem41} Assume that $\operatorname{dom}F_{u}\neq \emptyset $, for all $u\in
U,$ holds and let $x^{\ast }\in X^{\ast }$ be such that 
\begin{equation*}
q(x^{\ast })=\min\limits_{\begin{subarray}{l} u \in U \\ y_u^* \in Y_u^*
\end{subarray}}F_{u}^{\ast }(x^{\ast },y_{u}^{\ast }).
\end{equation*}%
Then there exist $u\in U$, $y_{u}^{\ast }\in Y_{u}^{\ast }$ such that 
\begin{equation*}
S^{\varepsilon }(x^{\ast })=B_{(u,y_{u}^{\ast })}^{\varepsilon }(x^{\ast
}),\ \forall \varepsilon \geq 0.
\end{equation*}
\end{lemma}

\noindent \textit{Proof.} By Lemma \ref{lem32} we have $B_{(u,y_{u}^{\ast
})}^{\varepsilon }(x^{\ast })\subset S^{\varepsilon }(x^{\ast })$.
Conversely, let $x\in S^{\varepsilon }(x^{\ast })$. By the exactness of $q$
at $x^{\ast }$, there exist $u\in U$ and $y_{u}^{\ast }\in Y_{u}^{\ast }$
such that 
\begin{equation*}
p(x)-\langle x^{\ast },x\rangle \leq -F_{u}^{\ast }(x^{\ast },y_{u}^{\ast
})+\varepsilon .
\end{equation*}%
Since $p(x)\in \mathbb{R}$ and $\operatorname{dom}F_{u}\neq \emptyset $, for all $%
u\in U,$ we have $F_{u}^{\ast }(x^{\ast },y_{u}^{\ast })\in \mathbb{R}$, $%
F_{u}(x,0_{u})\in \mathbb{R}$, 
\begin{equation*}
\Big(p(x)-F_{u}(x,0_{u})\Big)+\Big(F_{u}(x,0_{u})+F_{u}^{\ast }(x^{\ast
},y_{u}^{\ast })-\langle x^{\ast },x\rangle \Big)\leq \varepsilon .
\end{equation*}%
Consequently, there exist $\varepsilon _{1}\geq 0,\varepsilon _{2}\geq 0$
such that $\varepsilon _{1}+\varepsilon _{2}=\varepsilon ,$ 
\begin{equation*}
p(x)-F_{u}(x,0_{u})\leq \varepsilon _{1}\ \text{and\ }F_{u}(x,0_{u})+F_{u}^{%
\ast }(x^{\ast },y_{u}^{\ast })-\langle x^{\ast },x\rangle \leq \varepsilon
_{2},
\end{equation*}%
that is, $x\in J^{\varepsilon _{1}}(u)\ \text{and }(x,0_{u})\in
(M^{\varepsilon _{2}}F_{u})(x^{\ast },y_{u}^{\ast }) $. Thus, $x\in
A_{(u,y_{u}^{\ast })}^{(\varepsilon _{1},\varepsilon _{2})}(x^{\ast
})\subset B_{(u,y_{u}^{\ast })}^{\varepsilon }(x^{\ast })$, since $%
\varepsilon _{1}+\varepsilon _{2}=\varepsilon .$ \qed

\begin{theorem}[Strong robust duality]
\label{thm41} Assume that $\operatorname{dom}p\neq \emptyset $ and let $x^{\ast
}\in X^{\ast }$. The next statements are equivalent:

$\mathrm{(i)} $ \ $\inf {\mathrm{(RP)}_{{x^{\ast }}}}=\max {(\mathrm{ODP})_{{%
x^{\ast }}}}$,

$\mathrm{(ii)}$ $\exists u\in U,\ \exists y_{u}^{\ast }\in Y_{u}^{\ast }:$ $%
\left( M^{\varepsilon }p\right) (x^{\ast })=B_{(u,y_{u}^{\ast
})}^{\varepsilon }(x^{\ast }),\forall \varepsilon \geq 0$,

$\mathrm{(iii)} $ $\exists \bar{\varepsilon}>0,\ \exists u\in U,\exists
y_{u}^{\ast }\in Y_{u}^{\ast }:$ $\left( M^{\varepsilon }p\right) (x^{\ast
})=B_{(u,y_{u}^{\ast })}^{\varepsilon }(x^{\ast }),\forall \varepsilon \in
]0,\bar{\varepsilon}[$.
\end{theorem}

\textit{Proof.} Observe firstly that (i) means that 
\begin{equation*}
p^{\ast }(x^{\ast })=q(x^{\ast })=\min\limits_{\begin{subarray}{l} u \in U
\\ y_u^* \in Y_u^* \end{subarray}}F_{u}^{\ast }(x^{\ast },y_{u}^{\ast }).
\end{equation*}%
As $\operatorname{dom}p\neq \emptyset $, \eqref{32} holds, and then by Lemmas \ref%
{lem31} and \ref{lem41}, $\mathrm{(i)}$ implies the remaining conditions,
which are equivalent to each other, and also that $\mathrm{(iii)}$ implies $%
p^{\ast }(x^{\ast })=q(x^{\ast })$.

We now prove that $\mathrm{(iii)}$ implies $q(x^{\ast })=F_{u}^{\ast
}(x^{\ast },y_{u}^{\ast })$. Assume by contradiction that there exists $%
\varepsilon >0$ such that $q(x^{\ast })+\varepsilon <F_{u}^{\ast }(x^{\ast
},y_{u}^{\ast })$, and without loss of generality one can take $\varepsilon
\in \left] 0,\bar{\varepsilon}\right[ ,$ where $\bar{\varepsilon}>0$
appeared in (iii). Then, by $\mathrm{(iii)}$, $\left( M^{\varepsilon
}p\right) (x^{\ast })=B_{(u,y_{u}^{\ast })}^{\varepsilon }(x^{\ast }).$

Pick $x\in \left( M^{\varepsilon }p\right) (x^{\ast })=B_{(u,y_{u}^{\ast
})}^{\varepsilon }(x^{\ast }).$ Then, there are $\varepsilon _{1}\geq
0,\varepsilon _{2}\geq 0$, $\varepsilon _{1}+\varepsilon _{2}=\varepsilon $
and $x\in J^{\varepsilon _{1}}(u)$ and $(x,0_{u})\in (M^{\varepsilon
_{2}}F_{u})(x^{\ast },y_{u}^{\ast })$. In other words, 
\begin{eqnarray}
&&p(x)\leq F_{u}(x,0_{u})+\varepsilon _{1},  \label{eq41a} \\
&&F^{\ast }((x^{\ast },y_{u}^{\ast })+F_{u}(x,0_{u})\leq \langle x^{\ast
},x\rangle +\varepsilon _{2}.  \label{eq41b}
\end{eqnarray}%
It now follows from \eqref{eq41a}-\eqref{eq41b} that 
\begin{eqnarray*}
p^{\ast }\left( x^{\ast }\right) &\geq &\langle x^{\ast },x\rangle -p\left(
x\right) \geq\langle x^{\ast },x\rangle -F_{u}(x,0_{u})-\varepsilon _{1} \\
&\geq &\langle x^{\ast },x\rangle +F_{u}^{\ast }(x^{\ast },y_{u}^{\ast
})-\langle x^{\ast },x\rangle -\varepsilon _{2}-\varepsilon _{1}=F_{u}^{\ast
}(x^{\ast },y_{u}^{\ast })-\varepsilon >q(x^{\ast }),
\end{eqnarray*}%
which contradicts the fact that $p^{\ast }(x^{\ast })=q(x^{\ast })$. \qed

In deterministic optimization with linear perturbations we get the next
consequence from Theorem \ref{thm41}.

\begin{corollary}[Strong robust duality for Case 1]
\label{corol41} Let $F:X\times Y\rightarrow \mathbb{R}_{\infty }$, $%
p=F(\cdot ,0_{Y})$, and assume that $\operatorname{dom}p\neq \emptyset $. Then,
for each $x^{\ast }\in X^{\ast }$, {\ the strong duality for $\mathrm{(P)}%
_{x^{\ast }}$ in Case 1 holds at $x^{\ast }$, i.e., 
\begin{equation*}
\mathop {\inf }\limits_{x\in X}\left\{ {F\left( {x,{0_{Y}}}\right)
-\left\langle {{x^{\ast }},x}\right\rangle }\right\} =\mathop {\max }%
\limits_{{y^{\ast }}\in {Y^{\ast }}}-{F^{\ast }}\left( {{x^{\ast }},{y^{\ast
}}}\right) ,
\end{equation*}%
}if and only if one of the (equivalent) conditions (ii) or (iii) in Theorem %
\ref{thm41} holds with $B_{(u,y_{u}^{\ast })}^{\varepsilon }(x^{\ast })$
being replaced by 
\begin{equation}
B_{y^{\ast }}^{\varepsilon }(x^{\ast }):=\big\{x\in X\,:\,(x,0_{Y})\in
(M^{\varepsilon }F)(x^{\ast },y^{\ast })\big\}.  \label{By-star}
\end{equation}
\end{corollary}

\textit{Proof.} It is worth observing that we are in the non-uncertainty
case (i.e., $U$ is a singleton), and the set $B_{(u,y_{u}^{\ast
})}^{\varepsilon }(x^{\ast })$ writes as in \eqref{By-star} for each $%
(x^{\ast },y^{\ast })\in X^{\ast }\times Y^{\ast }$, $\varepsilon \geq 0$.
The conclusion follows from Theorem \ref{thm41}. \qed 

In the non-perturbation case, Theorem \ref{thm41} gives rise to

\begin{corollary}[Strong robust duality for Case 2]
\label{corol42} Let $(f_{u})_{u\in U}\subset \mathbb{R}_{\infty }^{X}$, $%
x^{\ast }\in X^{\ast }$, and $p=\sup\limits_{u\in U}f_{u}$ such that $%
\operatorname{dom}p\neq \emptyset $. Then, the robust duality formula 
\begin{equation*}
\Big(\sup\limits_{u\in U}f_{u}\Big)^{\ast }(x^{\ast })=\min\limits_{u\in
U}f_{u}^{\ast }(x^{\ast })
\end{equation*}%
holds if and only if one of the (equivalent) conditions (ii) or (iii) in
Theorem \ref{thm41} holds with $B_{(u,y_{u}^{\ast })}^{\varepsilon }(x^{\ast
})$ being replaced by 
\begin{equation}
B_{u}^{\varepsilon }(x^{\ast }):=\bigcup\limits_{\QATOP{\scriptstyle{%
\varepsilon _{1}}+{\varepsilon _{2}}=\varepsilon \hfill }{\scriptstyle{%
\varepsilon _{1}}\geqslant 0,{\varepsilon _{2}}\geqslant 0\hfill }}\Big(%
J^{\varepsilon _{1}}(u)\cap (M^{\varepsilon _{2}}f_{u})(x^{\ast })\Big).
\label{41}
\end{equation}
\end{corollary}

\textit{Proof.} Let $F_{u}(x,y_{u})=f_{u}(x)$, $p=\mathop{\sup }%
\limits_{u\in U}{f_{u}}$, and, from \eqref{34} and \eqref{35} (see the proof
of Corollary \ref{corol32}), 
\begin{equation*}
B_{(u,y_{u}^{\ast })}^{\varepsilon }({x^{\ast }})=\left\{ 
\begin{array}{ll}
\bigcup\limits_{\QATOP{\scriptstyle{\varepsilon _{1}}+{\varepsilon _{2}}%
=\varepsilon \hfill }{\scriptstyle{\varepsilon _{1}}\geqslant 0,{\varepsilon
_{2}}\geqslant 0\hfill }}\Big( J^{\varepsilon _{1}}(u)\cap (M^{\varepsilon
_{2}}f_{u})(x^{\ast })\Big) , & \mathrm{if}\ \ y_{u}^{\ast }=0_{u}^{\ast },
\\ 
\ \ \ \ \ \emptyset , & \text{else},%
\end{array}%
\right.
\end{equation*}
which in our situation, collapses to the set $B_{u}^{\varepsilon }(x^{\ast
}) $ defined by \eqref{41}. The conclusion now follows from Theorem \ref%
{thm41}. \qed

\section{Reverse strong and min-max robust duality}

Given $F_{u}:X\times Y_{u}\rightarrow (\mathbb{R}_{\infty })^{X}$ for each $%
u\in U$, $p=\sup\limits_{u\in U}F_{u}(\cdot ,0_{u})$, and $x^{\ast }\in
X^{\ast }$, we assume in this section that the problem $\mathrm{(RP)}%
_{x^{\ast }}$ is finite-valued and admits an optimal solution or, in other
words, that $\mathrm{argmin}(p-x^{\ast })=(M^{0}p)(x^{\ast })\neq \emptyset $%
. For convenience, we set 
\begin{equation*}
(Mp)(x^{\ast }):=(M^{0}p)(x^{\ast }),\ \ S(x^{\ast }):=S^{0}(x^{\ast }),\ 
\mathrm{and}
\end{equation*}%
\begin{equation}
\mathcal{A}(x^{\ast }):=\mathcal{A}^{0}(x^{\ast })=\bigcap_{\eta
>0}\bigcup\limits_{\QATOP{u\in U\hfill }{y_{u}^{\ast }\in Y_{u}^{\ast
}\hfill }}{{\bigcup\limits_{\QATOP{\scriptstyle{\varepsilon _{1}}+{%
\varepsilon _{2}}=\eta \hfill }{\scriptstyle{\varepsilon _{1}}\geq 0,\;{%
\varepsilon _{2}}\geq 0\hfill }}}}\Big\{x\in J^{\varepsilon
_{1}}(u)\,:\,(x,0_{u})\in (M^{\varepsilon _{2}}F_{u})(x^{\ast },y_{u}^{\ast
})\Big\}.  \label{51}
\end{equation}

\begin{theorem}[Reverse strong robust duality]
\label{thm51} Let $x^{\ast }\in X^{\ast }$ be such that $(Mp)(x^{\ast })\neq
\emptyset $ and let $\mathcal{A}(x^*)$ be as in \eqref{51}. The next
statements are equivalent:

$\mathrm{(i)}$ $\min {\mathrm{(RP)}_{{x^{\ast }}}}=\sup {(\mathrm{ODP})_{{%
x^{\ast }}}}$,

$\mathrm{(ii)} $ $(Mp)(x^{\ast }) = \mathcal{A} (x^{\ast })$.
\end{theorem}

\noindent \textit{Proof.} \ Since $(Mp)(x^{\ast })\neq \emptyset $, $%
\operatorname{dom}p\neq \emptyset $. It follows from Theorem \ref{thm31} that $%
\left[ {\mathrm{(i)}\Longrightarrow \mathrm{(ii)}}\right] $. For the
converse, let us pick $x\in (Mp)(x^{\ast })$. Then by (ii), for each $\eta
>0 $ there exist $u\in U$, $y_{u}^{\ast }\in Y_{u}^{\ast }$, $\varepsilon
_{1}\geq 0,\varepsilon _{2}\geq 0$ such that $\varepsilon _{1}+\varepsilon
_{2}=\eta $, $x\in J^{\varepsilon _{1}}(u)$, $(x,0_{u})\in (M^{\varepsilon
_{2}}F_{u})(x^{\ast },y_{u}^{\ast })$ and we have 
\begin{eqnarray*}
q(x^{\ast }) &\leq &F_{u}^{\ast }(x^{\ast },y_{u}^{\ast })\leq \langle
x^{\ast },x\rangle -F_{u}(x,0_{u})+\varepsilon _{2} \\
&\leq &\langle x^{\ast },x\rangle -p(x)+\varepsilon _{1}+\varepsilon
_{2}\leq p^{\ast }(x^{\ast })+\eta .
\end{eqnarray*}%
Since $\eta >0$ is arbitrary we get $q(x^{\ast })\leq p^{\ast }(x^{\ast })$,
which, together with the weak duality (see (\ref{WD})), yields $q(x^{\ast })
= \langle x^{\ast },x\rangle -p(x) = p^{\ast }(x^{\ast }) $, i.e., (i) holds
and we are done. \qed

In the deterministic case we obtain from Theorem \ref{thm51}:

\begin{corollary}[Reverse strong robust duality for Case 1]
\label{corol51} Let $F:X\times Y\rightarrow \mathbb{R}_{\infty }$, $x^{\ast
}\in X^{\ast }$, $p=F(\cdot ,0_{Y})$, and 
\begin{equation*}
\mathcal{A}(x^{\ast })=\bigcap_{\eta >0}\bigcup\limits_{y^{\ast }\in Y^{\ast
}}\Big\{x\in X\,:\,(x,0_{Y})\in (M^{\eta }F)(x^{\ast },y^{\ast })\Big\}.
\end{equation*}%
Assume that $(Mp)(x^{\ast })\neq \emptyset $. Then the next statements are
equivalent:

$\mathrm{(i)}$ $\mathop {\min }\limits_{x\in X}\left\{ {F\left( {x,{0_{Y}}}%
\right) -\left\langle {{x^{\ast }},x}\right\rangle }\right\} =\mathop {\sup }%
\limits_{{y^{\ast }}\in {Y^{\ast }}}-{F^{\ast }}\left( {{x^{\ast }},{y^{\ast
}}}\right)$,

$\mathrm{(ii)} $ $(Mp)(x^{\ast }) = \mathcal{A} (x^{\ast })$.
\end{corollary}

\begin{corollary}[Reverse strong robust duality for Case 2]
\label{corol53} Let $(f_{u})_{u\in U}\subset \mathbb{R}_{\infty }^{X}$, $%
p=\sup\limits_{u\in U}f_{u}$, $x^{\ast }\in X^{\ast }$, and 
\begin{equation*}
\mathcal{A}(x^{\ast }):=\bigcap_{\eta >0}\bigcup\limits_{u\in
U}\bigcup\limits_{\QATOP{\scriptstyle{\varepsilon _{1}}+{\varepsilon _{2}}%
=\eta \hfill }{\scriptstyle{\varepsilon _{1}}\geq 0,\;{\varepsilon _{2}}\geq
0\hfill }}\Big(J^{\varepsilon _{1}}(u)\cap (M^{\varepsilon
_{2}}f_{u})(x^{\ast })\Big),
\end{equation*}%
where 
\begin{equation*}
J^{\varepsilon _{1}}(u)=\Big\{x\in p^{-1}(\mathbb{R})\,:\,f_{u}(x)\geq
p(x)-\varepsilon _{1}\Big\}.
\end{equation*}%
Assume that $(Mp)(x^{\ast })\neq \emptyset $. Then the next statements are
equivalent:

$\mathrm{(i)}$ $\Big(\sup\limits_{u\in U}f_{u}\Big)^{\ast }(x^{\ast
})=\inf\limits_{u\in U}f_{u}^{\ast }(x^{\ast }),$\ with attainment at the
first member,

$\mathrm{(ii)} $ $(Mp)(x^{\ast }) = \mathcal{A} (x^{\ast })$.
\end{corollary}

Now, for each $u\in U$, $y_{u}^{\ast }\in Y_{u}^{\ast }$, $x^{\ast }\in
X^{\ast },$ we set%
\begin{equation*}
J(u):=J^{0}(u)=\Big\{x\in p^{-1}(\mathbb{R})\,:\,F_{u}(x,0_{u})=p(x)\Big\},
\end{equation*}%
\begin{equation*}
(MF_{u})(x^{\ast },y_{u}^{\ast }):=(M^{0}F_{u})(x^{\ast },y_{u}^{\ast })=%
\mathrm{argmin}\Big(F_{u}-\langle x^{\ast },\cdot \rangle -\langle
y_{u}^{\ast },\cdot \rangle \Big),
\end{equation*}%
and 
\begin{equation}
B_{(u,y_{u}^{\ast })}(x^{\ast }):=B_{(u,y_{u}^{\ast })}^{0}(x^{\ast })=\Big\{%
x\in J(u)\,:\,(x,0_{u})\in (MF_{u})(x^{\ast },y_{u}^{\ast })\Big\}.
\label{52}
\end{equation}

\begin{theorem}[Min-max robust duality]
\label{thm52} Let $x^{\ast }\in X^{\ast }$ be such that $(Mp)(x^{\ast })\neq
\emptyset $. The next statements are equivalent:

$\mathrm{(i)}$ {\ $\min \mathrm{(RP)}_{x^{\ast }}=\max \mathrm{(ODP)}%
_{x^{\ast }}$,}

$\mathrm{(ii)}$ $\exists u\in U,$ $\exists y_{u}^{\ast }\in Y_{u}^{\ast }:$ $%
\ (Mp)(x^{\ast })=B_{(u,y_{u}^{\ast })}(x^{\ast })$,

where $B_{(u,y_{u}^{\ast })}(x^{\ast })$ is the set defined in \eqref{52}.
\end{theorem}

\noindent \textit{Proof.} By Theorem \ref{thm41} we know that $[\mathrm{(i)}%
\Longrightarrow \mathrm{(ii)}]$. We now prove that $[\mathrm{(ii)}%
\Longrightarrow \mathrm{(i)}]$. Pick $x\in (Mp)(x^{\ast })$ which is
non-empty by assumption. Then by $\mathrm{(ii)}$, $x\in B_{(u,y_{u}^{\ast
})}(x^{\ast })$, which yields $x\in J(u)$ and $(x,0_{u})\in (MF_{u})(x^{\ast
},y_{u}^{\ast })$. Hence, 
\begin{eqnarray*}
q(x^{\ast }) &\leq &F_{u}^{\ast }(x^{\ast },y_{u}^{\ast })\leq \langle
x^{\ast },x\rangle -F_{u}(x,0_{u}) \\
&\leq &\langle x^{\ast },x\rangle -p(x)\leq p^{\ast }(x^{\ast })\leq
q(x^{\ast }),
\end{eqnarray*}%
which means that $q(x^{\ast })=F_{u}^{\ast }(x^{\ast },y_{u}^{\ast
})=\langle x^{\ast },x\rangle -p(x)=p^{\ast }(x^{\ast })$ and $\mathrm{(i)}$
follows. \qed

\begin{corollary}[Min-max robust duality for Case 1]
\label{corol52} Let $F:X\times Y\rightarrow \mathbb{R}_{\infty }$, $x^{\ast
}\in X^{\ast }$, $p=F(\cdot ,0_{Y})$, and for each $y^{\ast }\in Y^{\ast }$, 
\begin{equation*}
B_{y^{\ast }}(x^{\ast }):=\Big\{x\in X\,:\,(x,0_{Y})\in (MF)(x^{\ast
},y^{\ast })\Big\}.
\end{equation*}%
Assume that $(Mp)(x^{\ast })\neq \emptyset $. The next statements are
equivalent:

$\mathrm{(i)} $ {$\mathop {\min }\limits_{x\in X}\left\{ {F\left( {x,{0_{Y}}}%
\right) -\left\langle {{x^{\ast }},x}\right\rangle }\right\} = 
\mathop {\max
}\limits_{{y^{\ast }}\in {Y^{\ast }}}-{F^{\ast }}\left( {{x^{\ast }},{%
y^{\ast }}}\right)$, }

$\mathrm{(ii)} $ $\exists y^{\ast }\in Y^{\ast }$:\, $\ (Mp)(x^{\ast }) =
B_{y^{\ast }}(x^{\ast })$.
\end{corollary}

\begin{corollary}[Min-max robust duality for Case 2]
\label{corol54} Let $(f_{u})_{u\in U}\subset \mathbb{R}_{\infty }^{X}$, $%
p=\sup\limits_{u\in U}f_{u}$, $x^{\ast }\in X^{\ast }$, and and for each $%
u\in U$, 
\begin{equation*}
B_{u}(x^{\ast }):=J(u)\cap (Mf_{u})(x^{\ast }),
\end{equation*}%
where $J(u)=\{x\in p^{-1}(\mathbb{R})\,:\,f_{u}(x)=p(x)\}$. Assume that $%
(Mp)(x^{\ast })\neq \emptyset $. Then the next statements are equivalent:

$\mathrm{(i)} $ $\Big(\sup\limits_{u \in U} f_u\Big)^\ast (x^*)
=\min\limits_{u \in U} f_u^\ast (x^*)$, with attainment at the first member,

$\mathrm{(ii)} $ $\exists u \in U$:\, $\ (Mp)(x^{\ast }) = B_u(x^{\ast })$.
\end{corollary}

\section{Stable robust duality}

Let us first recall some notations. Given ${F_{u}}:X\times {Y_{u}}%
\rightarrow \mathbb{R}_{\infty }$, $u\in U,\ p=\mathop {\sup }\limits_{u\in
U}{F_{u}}(\cdot ,{0_{u}})$ and $q=\mathop {\inf }\limits_{\QATOP{%
\scriptstyle
u\in U\hfill }{\scriptstyle y_{u}^{\ast }\in Y_{u}^{\ast }\hfill }%
}F_{u}^{\ast }(\cdot ,y_{u}^{\ast })$. Remember that $p^{\ast }(x^{\ast
})\leq q(x^{\ast })$ for each $x^{\ast }\in X^{\ast }$. \textit{Stable
robust duality} means that $\inf {\mathrm{(RP)}_{{x^{\ast }}}}=\sup {(%
\mathrm{ODP})_{{x^{\ast }}}}$ for all $x^{\ast }\in X^{\ast }$, or
equivalently,

\begin{equation*}
p^{\ast }(x^{\ast })=q(x^{\ast }),\ \forall x^{\ast }\in X^{\ast }.
\end{equation*}%
Theorem \ref{thm31} says that, if $\operatorname{dom}p\neq \emptyset ,$ then
stable robust duality holds if and only if for each $\varepsilon \geq 0$ the
set-valued mappings $M^{\varepsilon }p,\ \mathcal{A}^\varepsilon:X^{\ast
}\rightrightarrows X$ coincide, where, for each $x^{\ast }\in X^{\ast }$, 
\begin{eqnarray*}
(M^{\varepsilon }p)(x^{\ast }):= &&\varepsilon -\mathrm{argmin}(p-x^{\ast }),
\\
\mathcal{A}^{\varepsilon }(x^{\ast }):= &&\bigcap_{\eta >0}\bigcup\limits_{%
\QATOP{u\in U\hfill }{y_{u}^{\ast }\in Y_{u}^{\ast }\hfill }}{{%
\bigcup\limits_{\QATOP{\scriptstyle{\varepsilon _{1}}+{\varepsilon _{2}}%
=\varepsilon +\eta \hfill }{\scriptstyle{\varepsilon _{1}}\geq 0,\;{%
\varepsilon _{2}}\geq 0\hfill }}}}\Big\{x\in J^{\varepsilon
_{1}}(u)\,:\,(x,0_{u})\in (M^{\varepsilon _{2}}F_{u})(x^{\ast },y_{u}^{\ast
})\Big\}.
\end{eqnarray*}%
Consequently, stable robust duality holds if and only if for each $%
\varepsilon \geq 0$, the inverse set-valued mappings 
\begin{equation*}
(M^{\varepsilon }p)^{-1},\ (\mathcal{A}^{\varepsilon
})^{-1}:X\rightrightarrows X^{\ast },
\end{equation*}%
coincide. Recall that $(M^{\varepsilon }p)^{-1}$ is nothing but the $%
\varepsilon $-subdifferential of $p$ at $x$.

Let us now make explicit $(\mathcal{A}^{\varepsilon })^{-1}$. To this end we
need to introduce for each $\varepsilon \geq 0$ the ($\varepsilon $-active
indexes) set-valued mapping $I^{\varepsilon }:X\rightrightarrows U$ with 
\begin{equation}  \label{I-ep-x}
\ I^{\varepsilon }(x)=\left\{ 
\begin{array}{ll}
\Big\{u\in U\,:\,F_{u}(x,0_{u})\geq p(x)-\varepsilon \Big\}, & \mathrm{if}%
\quad p(x)\in \mathbb{R}, \\ 
\ \ \ \emptyset , & \mathrm{if}\quad p(x)\not\in \mathbb{R}.%
\end{array}%
\right.
\end{equation}%
We observe that $I^{\varepsilon }$ is nothing but the inverse of the
set-valued mapping $J^{\varepsilon }:U\rightrightarrows X$ defined in %
\eqref{J-ep}.

\begin{lemma}
\label{lem61} For each $(\varepsilon ,x)\in \mathbb{R}_{+}\times X$ one has 
\begin{equation*}
{(\mathcal{A}^{\varepsilon })^{-1}}(x)=\bigcap\limits_{\eta >0}{%
\bigcup\limits_{\QATOP{\scriptstyle{\varepsilon _{1}}+{\varepsilon _{2}}%
=\varepsilon +\eta \hfill }{\scriptstyle{\varepsilon _{1}}\geqslant 0,{%
\varepsilon _{1}}\geqslant 0\hfill }}{\bigcup\limits_{u\in {I^{\varepsilon
_{1}}}(x)}\mathrm{proj}_{X^{\ast }}^{u}{\partial ^{{\varepsilon _{2}}}}{F_{u}%
}(x,{0_{u}})}},
\end{equation*}%
where ${{\mathrm{proj}_{X^{\ast }}^{u}:}}X^{\ast }\times Y_{u}^{\ast
}\longrightarrow X^{\ast }$ is the projection mapping $\mathrm{proj}%
_{X^{\ast }}^{u}(x^{\ast },y_{u}^{\ast })=x^{\ast }.$
\end{lemma}

\noindent \textit{Proof.} Let $(\varepsilon ,x,x^{\ast })\in \mathbb{R}%
_{+}\times X\times X^{\ast }$. One has 
\begin{eqnarray*}
x^{\ast }\in (\mathcal{A}^{\varepsilon })^{-1}(x) &\Leftrightarrow & x\in 
\mathcal{A}^{\varepsilon }(x^{\ast }) \\
&\Leftrightarrow &\left\{ 
\begin{array}{l}
\forall \eta >0,\exists u\in U,\exists y_{u}^{\ast }\in Y_{u}^{\ast
},\exists (\varepsilon _{1},\varepsilon _{2})\in \mathbb{R}_{+}^{2}\text{ }%
\mathrm{such\ that} \\ 
\varepsilon _{1}+\varepsilon _{2}=\varepsilon +\eta ,\text{ }x\in
J^{\varepsilon _{1}}(u)\ \mathrm{and}\ (x,0_{u})\in (M^{\varepsilon
_{2}}F_{u})(x^{\ast },y_{u}^{\ast })%
\end{array}%
\right. \\
&\Leftrightarrow & \left\{ 
\begin{array}{l}
\forall \eta >0,\exists u\in U,\exists y_{u}^{\ast }\in Y_{u}^{\ast
},\exists (\varepsilon _{1},\varepsilon _{2})\in \mathbb{R}_{+}^{2}\text{ }%
\mathrm{such\ that\ } \\ 
\varepsilon _{1}+\varepsilon _{2}=\varepsilon +\eta ,\text{ }u\in
I^{\varepsilon _{1}}(x),\ \mathrm{and}\ (x^{\ast },y_{u}^{\ast })\in \Big(%
\partial ^{\varepsilon _{2}}F_{u}\Big)(x,0_{u})%
\end{array}%
\right. \\
&\Leftrightarrow &\left\{ 
\begin{array}{l}
\forall \eta >0,\exists u\in U,\exists (\varepsilon _{1},\varepsilon
_{2})\in \mathbb{R}_{+}^{2}\text{ }\mathrm{\ such\ that\ } \\ 
\mathrm{\ }\varepsilon _{1}+\varepsilon _{2}=\varepsilon +\eta ,\text{ }u\in
I^{\varepsilon _{1}}(x),\ \mathrm{and}\ x^{\ast }\in \mathrm{proj}_{X^{\ast
}}^{u}\Big(\partial ^{\varepsilon _{2}}F_{u}\Big)(x,0_{u})%
\end{array}%
\right. \\
& \Leftrightarrow & x^{\ast} \in \bigcap\limits_{\eta>0} \bigcup\limits_{\QATOP{\scriptstyle{\varepsilon_{1}}+{\varepsilon _{2}}=\varepsilon +\eta
\hfill }{\scriptstyle{\varepsilon {1}}\geqslant 0, \varepsilon_{1} \geqslant 0\hfill}} \bigcup\limits_{u\in {I^{\varepsilon_{1}}}(x)} \mathrm{proj}_{X^{\ast}}^{u}(\partial^{\varepsilon_{2}} F_{u})(x,0_{u}).
\end{eqnarray*}
\qed

\medskip

Now, for each $(\varepsilon, x) \in \mathbb{R}_+ \times X$, let us set 
\begin{equation}  \label{61}
C^\varepsilon (x) := \bigcap\limits_{\eta >0} \bigcup\limits_{\QATOP{%
\scriptstyle{\varepsilon _{1}}+{\varepsilon _{2}}=\varepsilon +\eta \hfill }{%
\scriptstyle{\varepsilon _{1}}\geqslant 0,{\varepsilon _{1}}\geqslant
0\hfill }}\bigcup\limits_{u\in I^{\varepsilon _{1}}(x)}\mathrm{proj}^u
_{X^{\ast }}(\partial ^{\varepsilon _{2}}F_{u})(x,{0_{u}}).
\end{equation}

Applying Theorem \ref{thm31} and Lemma \ref{lem61} we obtain:

\begin{theorem}[Stable robust duality]
\label{thm61} Assume that $\operatorname{dom}p\neq \emptyset $. The next
statements are equivalent:

$\mathrm{(i)} $ $\inf {\mathrm{(RP)}_{{x^{\ast }}}}=\sup {(\mathrm{ODP})_{{%
x^{\ast }}}}$ for all $x^* \in X^\ast$,

$\mathrm{(ii)} $ $\partial^\varepsilon p(x) = C^\varepsilon (x), \ \forall
(\varepsilon, x) \in \mathbb{R}_+ \times X$,

$\mathrm{(iii)} $ $\exists \bar\varepsilon > 0$: \ $\partial^\varepsilon
p(x) = C^\varepsilon (x), \ \forall (\varepsilon, x) \in ]0,
\bar\varepsilon[ \times X$.
\end{theorem}

\begin{corollary}[Stable robust duality for Case 1]
\label{corol61} Let $F:X\times Y\rightarrow \mathbb{R}_{\infty }$ be such
that $\operatorname{dom}F(\cdot ,0_{Y})\neq \emptyset $. Let ${{\mathrm{proj}%
_{X^{\ast }}:}}X^{\ast }\times Y^{\ast }\longrightarrow X^{\ast }$ be the
projection mapping $\mathrm{proj}_{X^{\ast }}(x^{\ast },y^{\ast })=x^{\ast
}. $ Then, the next statements are equivalent:

$\mathrm{(i)} $ $\inf\limits_{x \in X} \Big\{ F(x, 0_Y) - \langle x^*,
x\rangle \Big\} = \sup\limits_{y^* \in Y^\ast} - F^\ast (x^*, y^*), \
\forall x^* \in X^\ast$,

$\mathrm{(ii)}$ $(\partial ^{\varepsilon }p)(x)=\bigcap_{\eta >0}\mathrm{proj%
}_{X^{\ast }}(\partial ^{{\varepsilon }}{F})(x,{0_{Y}}),\ \forall
(\varepsilon ,x)\in \mathbb{R}_{+}\times X$,

$\mathrm{(iii)}$ $\exists \bar{\varepsilon}>0$: \ $(\partial ^{\varepsilon
}p)(x)=\bigcap_{\eta >0}\mathrm{proj}_{X^{\ast }}(\partial ^{{\varepsilon }}{%
F})(x,{0_{Y}}),\ \forall (\varepsilon ,x)\in ]0,\bar{\varepsilon}[\times X$.
\end{corollary}

\textit{Proof.} Let $U=\{u_{0}\}$ and $F=F_{u_{0}}:X\times Y\rightarrow 
\mathbb{R}_{\infty }$, $Y=Y_{u_{0}}$, $p=F(\cdot ,0_{Y})$. Then for each $%
(\varepsilon ,x)\in \mathbb{R}_{+}\times X$, 
\begin{equation*}
I^{\varepsilon }(x)=\left\{ 
\begin{array}{ll}
\{u_{0}\},\ \  & \mathrm{if}\quad p(x)\in \mathbb{R}, \\ 
\emptyset , & \mathrm{if}\quad p(x)\not\in \mathbb{R},%
\end{array}%
\right.
\end{equation*}%
and%
\begin{equation}
\bigcup\limits_{\QATOP{\scriptstyle{\varepsilon _{1}}+{\varepsilon _{2}}%
=\varepsilon \hfill }{\scriptstyle{\varepsilon _{1}}\geqslant 0,{\varepsilon
_{1}}\geqslant 0\hfill }}\bigcup\limits_{u\in I^{\varepsilon _{1}}(x)}%
\mathrm{proj}_{X^{\ast }}^{u}(\partial ^{\varepsilon _{2}}F_{u})(x,{0_{u}})=%
\mathrm{proj}_{X^{\ast }}\Big(\partial ^{\varepsilon }F\Big)(x,0_{Y}).
\label{62}
\end{equation}%
The conclusion now follows from \eqref{61}-\eqref{62} and Theorem \ref{thm61}%
. \qed

\begin{remark}
Condition $\mathrm{(ii)}$ in Corollary \ref{corol61} was quoted in \cite[%
Theorem 4.3]{Grad16} for all $(\varepsilon ,x)\in ]0,+\infty \lbrack \times 
\mathbb{R},$ which is equivalent.
\end{remark}

\begin{corollary}[Stable robust duality for Case 2]
\label{corol62} Let $(f_{u})_{u\in U}\subset \mathbb{R}_{\infty }^{X},$ $%
p=\sup\limits_{u\in U}f_{u}$ and assume that $\operatorname{dom}p\neq \emptyset $%
. The next statements are equivalent:

$\mathrm{(i)} $ $\Big(\sup\limits_{u \in U} f_u\Big)^\ast (x^*)
=\inf\limits_{u \in U} f_u^\ast (x^*)$, \ $\forall x^* \in X^\ast$,

$\mathrm{(ii)} $ $(\partial^\varepsilon p)(x) = C^\varepsilon (x), \ \forall
(\varepsilon, x) \in \mathbb{R}_+ \times X$,

$\mathrm{(iii)}$ $\exists \bar{\varepsilon}>0$: \ $(\partial ^{\varepsilon
}p)(x)=C^{\varepsilon }(x),\ \forall (\varepsilon ,x)\in ]0,\bar{\varepsilon}%
[\times X$, \newline
where $C^{\varepsilon }(x)$ is the set 
\begin{equation}
C^{\varepsilon }(x)=\bigcap\limits_{\eta >0}\bigcup\limits_{\QATOP{%
\scriptstyle{\varepsilon _{1}}+{\varepsilon _{2}}=\varepsilon +\eta \hfill }{%
\scriptstyle{\varepsilon _{1}}\geqslant 0,{\varepsilon _{1}}\geqslant
0\hfill }}\bigcup\limits_{u\in I^{\varepsilon _{1}}(x)}(\partial
^{\varepsilon _{2}}f_{u})(x),\ \ \forall (\varepsilon ,x)\in \mathbb{R}%
_{+}\times X.  \label{63}
\end{equation}
\end{corollary}

\textit{Proof.} Let $F_{u}:X\times Y_{u}\rightarrow \mathbb{R}_{\infty }\ $%
be such that $F_{u}(x,y_{u})=f_{u}(x)$ for all $u\in U$. Then for any $%
(\varepsilon ,x)\in \mathbb{R}_{+}\times X$, 
\begin{equation*}
I^{\varepsilon }(x)=\left\{ 
\begin{array}{ll}
\Big\{u\in U\,:\,f_{u}(x)\geq p(x)-\varepsilon \Big\}, & \mathrm{if}\quad
p(x)\in \mathbb{R}, \\ 
\ \ \emptyset , & \mathrm{if}\quad p(x)\not\in \mathbb{R},%
\end{array}%
\right.
\end{equation*}%
\begin{equation*}
(\partial ^{\varepsilon }F_{u})(x,0_{u})=(\partial ^{\varepsilon
}f_{u})(x)\times \{0_{u}^{\ast }\},\ \ \ \forall (u,\varepsilon ,x)\in
U\times \mathbb{R}_{+}\times X,
\end{equation*}%
and $C^{\varepsilon }(x)$ reads as in \eqref{63}. The conclusion now follows
from Theorem \ref{thm61}. \qed

\section{Stable strong robust duality}

We retain all the notations used in the Sections 3-6. Given $(\varepsilon,
u) \in \mathbb{R}_+ \times U$ and $y_u^* \in Y_U^\ast$ we have introduced in
Section 4 the set-valued mapping $B_{(u,y_{u}^{\ast })}^{\varepsilon } :
X^\ast \rightrightarrows X$ defined by 
\begin{equation*}
B_{(u,y_{u}^{\ast })}^{\varepsilon }({x^{\ast }}) =\bigcup\limits_{\QATOP{%
\scriptstyle{\varepsilon _{1}}+{\varepsilon _{2}}=\varepsilon \hfill }{%
\scriptstyle{\varepsilon _{1}}\geqslant 0,{\varepsilon _{2}}\geqslant
0\hfill }}\Big\{x\in J^{\varepsilon _1}(u)\, :\, (x,0_{u}) \in ( M^{
\varepsilon _{2}} F_{u}) (x^\ast ,y_u^*)\Big\}.
\end{equation*}%
Let us now define $B^\varepsilon : X^\ast \rightrightarrows X$ by setting 
\begin{equation*}
B^\varepsilon (x^*) := \bigcup\limits_{\QATOP{u \in U}{y_u^* \in Y_u^\ast}}
B_{(u,y_{u}^{\ast })}^{\varepsilon }(x^*) , \ \ \forall x^* \in X^\ast.
\end{equation*}

\begin{lemma}
\label{lem71} For each $(\varepsilon, x) \in \mathbb{R}_+ \times X$ we have 
\begin{equation*}
(B^\varepsilon)^{-1} (x) = \bigcup\limits_{\QATOP{\scriptstyle{\varepsilon
_{1}}+{\varepsilon _{2}}=\varepsilon \hfill }{\scriptstyle{\varepsilon _{1}}%
\geqslant 0,{\varepsilon _{2}}\geqslant 0\hfill }} \bigcup\limits_{u\in
I^{\varepsilon _{1}}(x)}\mathrm{proj}^u_{X^{\ast }}(\partial ^{\varepsilon
_{2}}F_{u})(x,{0_{u}}).
\end{equation*}
\end{lemma}

\noindent \textit{Proof.} $x^* \in (B^\varepsilon)^{-1} (x)$ means that
there exist $u \in U$, $y_u^* \in Y_u^\ast$ $\varepsilon_1 \geq 0$, $%
\varepsilon_2 \geq 0$, such that $\varepsilon_1 + \varepsilon_2 =
\varepsilon $, $x \in J^{\varepsilon_1} (u) $, and $(x, 0_u) \in
(M^{\varepsilon_2}F_u)(x^*, y_u^*)$, or, equivalently, $u \in
I^{\varepsilon_1} (x)$, and $(x^*, y_u^*) \in
(\partial^{\varepsilon_2}F_u)(x, 0_u)$. In other words, there exist $u \in U$%
, $y_u^* \in Y_u^\ast$ such that $x \in B^\varepsilon_{(u, y_u^*)} (x^*)$,
that is $x \in B^\varepsilon (x^*)$. \qed 

For each $\varepsilon \geq 0$ let us introduce the set-valued mapping $%
D^{\varepsilon }:=(B^{\varepsilon })^{-1}$. Now Lemma \ref{lem71} writes: 
\begin{equation}
D^{\varepsilon }(x)\ =\ \bigcup\limits_{\QATOP{\scriptstyle{\varepsilon _{1}}%
+{\varepsilon _{2}}=\varepsilon \hfill }{\scriptstyle{\varepsilon _{1}}%
\geqslant 0,{\varepsilon _{2}}\geqslant 0\hfill }}\bigcup\limits_{u\in
I^{\varepsilon _{1}}(x)}\mathrm{proj}_{X^{\ast }}^{u}(\partial ^{\varepsilon
_{2}}F_{u})(x,{0_{u}}),\ \forall (\varepsilon ,x)\in \mathbb{R}_{+}\times X.
\label{71}
\end{equation}%
Note that 
\begin{equation}
C^{\varepsilon }(x)=\bigcap\limits_{\eta >0}D^{\varepsilon +\eta }(x),\
\forall (\varepsilon ,x)\in \mathbb{R}_{+}\times X,  \label{eq72}
\end{equation}%
and that $D^{\varepsilon }(x)=\emptyset $ whenever $p(x)\not\in \mathbb{R}$.

We now provide a characterization of stable strong robust duality in terms
of $\varepsilon$-subdifferential formulas.

\begin{theorem}[Stable strong robust duality]
\label{thm71} Assume that $\operatorname{dom} p \ne \emptyset$, and let $%
D^\varepsilon$ as in \eqref{71}. The next statements are equivalent:

$\mathrm{(i)} $ $\inf {\mathrm{(RP)}_{{x^{\ast }}}}=\max {(\mathrm{ODP})_{{%
x^{\ast }}}} = \max\limits_{\QATOP{u \in U\hfill }{y_u^* \in Y_u^\ast\hfill }%
} - F_u^\ast (x^*, y_u^*),\ \ \forall x^* \in X^\ast$,

{$\mathrm{(ii)} $} $\partial^\varepsilon p(x) = D^\varepsilon (x), \ \forall
(\varepsilon, x) \in \mathbb{R}_+ \times X$.
\end{theorem}

\noindent \textit{Proof.} $[\mathrm{(i)} \Longrightarrow \mathrm{(ii)}]$ Let 
$x^* \in \partial^\varepsilon p(x)$. Then $x \in (M^\varepsilon p)(x^*)$.
Since strong robust duality holds at $x^*$, Theorem \ref{thm41} says that
there exist $u \in U$, $y_u^* \in Y_u^\ast$ such that $x \in
B^\varepsilon_{(u, y_u^*)}(x^*) \subset B^\varepsilon(x^*)$, and finally $%
x^* \in D^\varepsilon (x)$ by Lemma \ref{lem71}. Thus $\partial^\varepsilon
p(x) \subset D^\varepsilon (x)$.

{\ Now, let $x^{\ast }\in D^{\varepsilon }(x)$. By Lemma \ref{lem71} we have 
$x\in B^{\varepsilon }(x^{\ast })$ and there exist $u\in U$, $y_{u}^{\ast
}\in Y_{u}^{\ast }$ such that $x\in B_{(u,y_{u}^{\ast })}^{\varepsilon
}(x^{\ast })$. By Lemma \ref{lem32} and the definition of $B_{(u,y_{u}^{\ast
})}^{\varepsilon }(x^{\ast })$ we have $x\in S^{\varepsilon }(x^{\ast })$,
and, by \eqref{31}, $x\in (M^{\varepsilon }p)(x^{\ast })$ which means that $%
x^{\ast }\in \partial ^{\varepsilon }p(x)$, and hence, $D^{\varepsilon
}(x)\subset \partial ^{\varepsilon }p(x)$. Thus (ii) follows.}

{\ $[\mathrm{(ii)}\Longrightarrow \mathrm{(i)}]$} If $p^{\ast }(x^{\ast
})=+\infty $ then $q(x^{\ast })=+\infty $ and one has $p^{\ast }(x^{\ast
})=F_{u}^{\ast }(x^{\ast },y_{u}^{\ast })=+\infty $ for all $u\in U$, $%
y_{u}^{\ast }\in Y_{u}^{\ast }$, and (i) holds. Assume that $p^{\ast
}(x^{\ast })\in \mathbb{R}$ and pick $x\in p^{-1}(\mathbb{R})$ which is
non-empty as $\operatorname{dom}p\neq \emptyset $ and $p^{\ast }(x^{\ast })\in 
\mathbb{R}$. Let $\varepsilon :=p(x)+p^{\ast }(x^{\ast })-\langle x^{\ast
},x\rangle $. Then $\varepsilon \geq 0$ and we have $x^{\ast }\in \partial
^{\varepsilon }p(x)$. By (ii) $x\in D^{\varepsilon }(x)$ and hence, there
exist $\varepsilon _{1}\geq 0,\varepsilon _{2}\geq 0$, $u\in U$, and $%
y_{u}^{\ast }\in Y_{u}^{\ast }$ such that $\varepsilon _{1}+\varepsilon
_{2}=\varepsilon $, $u\in I^{\varepsilon _{1}}(x)$, $(x^{\ast },y_{u}^{\ast
})\in (\partial ^{\varepsilon _{2}}F_{u})(x,0_{u})$. We have 
\begin{eqnarray*}
q(x^{\ast }) &\leq &F_{u}^{\ast }(x^{\ast },y_{u}^{\ast })\leq \langle
x^{\ast },x\rangle -F_{u}(x,0_{u})+\varepsilon _{2} \\
&\leq &\langle x^{\ast },x\rangle -p(x)+\varepsilon _{1}+\varepsilon
_{2}=p^{\ast }(x^{\ast })\ \ \ \ \ \ \mathrm{(by\ definition\ of\ }%
\varepsilon ) \\
&\leq &q(x^{\ast }),
\end{eqnarray*}%
and finally, $q(x^{\ast })=F_{u}^{\ast }(x^{\ast },y_{u}^{\ast })=p^{\ast
}(x^{\ast })$, which is (i).\qed

Next, as usual, we give two consequences of Theorem \ref{thm71} for the
non-uncertainty and non-parametric cases.

\begin{corollary}[Stable strong duality for Case 1]
\label{corol71} Let $F:X\times Y\rightarrow \mathbb{R}_{\infty }$, $%
p=F(\cdot ,0_{Y})$, $\operatorname{dom}p\neq \emptyset $. The next statements are
equivalent:

$\mathrm{(i)} $ $\inf\limits_{x \in X} \Big\{ F(x, 0_Y) - \langle x^*,
x\rangle\Big\}  = \max\limits_{ y^* \in Y^ \ast} - F^\ast (x^*, y^*), \
\forall x^* \in X^\ast$,

$\mathrm{(ii)}$\ \ $\partial ^{\varepsilon }p(x)=\mathrm{proj}_{X^{\ast
}}(\partial ^{\varepsilon }F)(x,{0_{y}}),\ \forall (\varepsilon ,x)\in 
\mathbb{R}_{+}\times X.$
\end{corollary}

\textit{Proof.} This is the non-uncertainty case (i.e., the uncertainty set
is a singleton) of the general problem $\mathrm{(RP)}_{x^{\ast }}$, with $%
U=\{u_{0}\}$ and $F_{u_{0}}=F:X\times Y\rightarrow \mathbb{R}_{\infty }$. We
have from \eqref{62}, 
\begin{equation}
D^{\varepsilon }(x)=\mathrm{proj}_{X^{\ast }}(\partial ^{\varepsilon }F)(x,{%
0_{Y}}),\ \forall (\varepsilon ,x)\in \mathbb{R}_{+}\times X.  \label{72}
\end{equation}%
%
%
%
%
%
%
%
%
%
%
%
%
%
%
%
%
%
%
%
%
%
%
%
%
%
%
%
%
%
The conclusion now follows from Theorem \ref{thm71}. \qed

\begin{corollary}[Stable strong duality for Case 2]
\label{corol72} Let $(f_{u})_{u\in U}\subset \mathbb{R}_{\infty }^{X}$, $%
p=\sup\limits_{u\in U}f_{u}$, and $\operatorname{dom}p\neq \emptyset $. The next
statements are equivalent:

$\mathrm{(i)} $ $(\sup\limits_{u \in U} f_u)^\ast(x^*) = \min\limits_{ u\in
U} f_u^\ast (x^*), \ \forall x^* \in X^\ast$,

$\mathrm{(ii)}$\ \ $\partial ^{\varepsilon }p(x)=D^{\varepsilon }(x),\
\forall (\varepsilon ,x)\in \mathbb{R}_{+}\times X$, where 
\begin{equation}
D^{\varepsilon }(x)=\bigcup\limits_{\QATOP{\scriptstyle{\varepsilon _{1}}+{%
\varepsilon _{2}}=\varepsilon \hfill }{\scriptstyle{\varepsilon _{1}}%
\geqslant 0,{\varepsilon _{2}}\geqslant 0\hfill }}\bigcup\limits_{u\in
I^{\varepsilon _{1}}(x)}(\partial ^{\varepsilon _{2}}f_{u})(x),\ \forall
(\varepsilon ,x)\in \mathbb{R}_{+}\times X,  \label{73}
\end{equation}%
and 
\begin{equation*}
I^{\varepsilon }(x)=\left\{ 
\begin{array}{ll}
\Big\{u\in U\,:\,f_{u}(x)\geq p(x)-\varepsilon \Big\} & \mathrm{if}\quad
p(x)\in \mathbb{R}, \\ 
\ \ \emptyset & \mathrm{if}\quad p(x)\not\in \mathbb{R}\text{.}%
\end{array}%
\right.
\end{equation*}
\end{corollary}

\textit{Proof.} In this non-parametric situation, let $%
F_{u}(x,y_{u})=f_{u}(x)$. It is easy to see that in this case, the set $%
D^{\varepsilon }(x)$ can be expressed as in\ \eqref{73}, and the conclusion
follows from Theorem \ref{thm71}. \qed

\section{Exact subdifferential formulas: Robust Basic Qualification condition%
}

Given ${F_{u}}:X\times {Y_{u}}\rightarrow \mathbb{R}_{\infty },\ u\in U$, as
usual, we let $p=\mathop {\sup }\limits_{u\in U}{F_{u}}(\cdot ,{0_{u}})$, $%
q:=\mathop {\inf }\limits_{\left( {u,y_{u}^{\ast }}\right) \in \Delta
}\;F_{u}^{\ast }(\cdot ,y_{u}^{\ast })$. Again, we consider the robust
problem $\mathrm{(RP)}_{x^{\ast }}$ and its robust dual problem $\mathrm{%
(ODP)}_{x^{\ast }}$ given in \eqref{RPx-star} and \eqref{ODPx-star},
respectively. Note that the reverse strong robust duality holds at $x^{\ast
} $ means that, for some $\bar{x}\in X$, it holds: 
\begin{equation}
-p^{\ast }(x^{\ast })=\min {\mathrm{(RP)}_{{x^{\ast }}}}=\sup_{u\in U}F_{u}(%
\bar{x},0_{u})-\langle x^{\ast },\bar{x}\rangle =p(\bar{x})-\langle x^{\ast
},\bar{x}\rangle =\sup {(\mathrm{ODP})_{{x^{\ast }}}}=-q(x^{\ast }).
\label{strongrobustduality}
\end{equation}

Now, let us set, for each $x\in X$, 
\begin{eqnarray}
D(x):= &&D^{0}(x)=\bigcup\limits_{u\in I(x)}\mathrm{proj}^u_{X^{\ast
}}(\partial F_{u})(x,{0_{u}}),  \label{81} \\
C(x):= &&C^{0}(x)=\bigcap\limits_{\eta >0}\bigcup\limits_{\QATOP{\scriptstyle%
{\varepsilon _{1}}+{\varepsilon _{2}}=\eta \hfill }{\scriptstyle{\varepsilon
_{1}}\geqslant 0,{\varepsilon _{2}}\geqslant 0\hfill }}\bigcup\limits_{u\in
I^{\varepsilon _{1}}(x)}\mathrm{proj}^u_{X^{\ast }}(\partial ^{\varepsilon
_{2}}F_{u})(x,{0_{u}}),  \label{82}
\end{eqnarray}%
where $I^{\varepsilon _{1}}(x)$ is defined as in \eqref{I-ep-x} and 
\begin{equation}
I(x):=\left\{ 
\begin{array}{ll}
\left\{ u\in U:F_{u}(x,0_{u})=p(x)\right\} , & \mathrm{if}\quad p(x)\in 
\mathbb{R}, \\ 
\emptyset , & \mathrm{if}\quad p(x)\not\in \mathbb{R}.%
\end{array}%
\right.  \label{I-x}
\end{equation}

\begin{lemma}
\label{lem81} For each $x \in X$, it holds 
\begin{equation*}
D(x) \subset C(x) \subset \partial p (x).
\end{equation*}
\end{lemma}

\noindent \textit{Proof.} The first inclusion is easy to check. Now let $x^*
\in C(x)$. For each $\eta > 0$ there exist $(\varepsilon_1, \varepsilon_2 )
\in \mathbb{R}_+^2$, $u \in I^{\varepsilon_1}(x)$, and $y_u^* \in Y_u^\ast$
such that $\varepsilon_1 + \varepsilon_2 = \eta$ and $(x^*, y_u^*) \in
(\partial^{\varepsilon_2} F_u)(x, 0_u)$. We then have $F_u^\ast (x^*, y_u^*)
+ F_u(x, 0_u) - \langle x^*, x\rangle \leq \varepsilon_2$, $p(x) \leq F_u(x,
0_u) + \varepsilon_1$ (as $u \in I^{\varepsilon_1}(x)$), and $p^\ast (x^*)
\leq q(x^*) \leq F^\ast_u(x^*, y^*_u)$. Consequently, 
\begin{equation*}
p^\ast (x^*) + p(x) - \langle x^*, x\rangle \leq F_u^\ast (x^*, y_u^*) +
F_u(x, 0_u) + \varepsilon_1- \langle x^*, x\rangle \leq\varepsilon_1 +
\varepsilon_2 = \eta .
\end{equation*}
Since $\eta > 0$ is arbitrary we get $p^\ast (x^*) + p(x) -\langle x^*,
x\rangle \leq 0$, which means that $x^* \in \partial p (x)$. The proof is
complete. \qed

\begin{theorem}
\label{thm81} Let $x \in p^{-1}(\mathbb{R})$ and $C(x)$ be as in \eqref{82}.
The next statements are equivalent:

$\mathrm{(i)} $\ \ $\partial p (x) = C(x)$,

$\mathrm{(ii)} $\ Reverse strong robust duality holds at each $x^* \in
\partial p (x) $,

$\mathrm{(iii)} $ Robust duality holds at each $x^* \in \partial p (x)$.
\end{theorem}

\noindent \textit{Proof.} $[ \mathrm{(i)} \Longrightarrow \mathrm{(ii)}]$
Let $x^* \in \partial p (x)$. We have $x^* \in C(x) = (\mathcal{A})^{-1} (x)$
(see Lemma \ref{lem61} with $\varepsilon = 0$). Then $x \in \mathcal{A}
(x^*)= S(x^*)$ (see \eqref{eqlem34} with $\varepsilon = 0$), and therefore, 
\begin{equation*}
- p^\ast (x^*) \leq p(x) - \langle x^*, x\rangle \leq - q (x^*) \leq -
p^\ast(x^*),
\end{equation*}
\begin{equation*}
- p^\ast (x^*) = \min\limits_{z \in X} \{ p(z) - \langle x^*, z\rangle \} =
p(x) - \langle x^*, x\rangle = -q(x^*) ,
\end{equation*}
that means that reverse strong robust duality holds at $x^*$ (see %
\eqref{strongrobustduality}).

$[ \mathrm{(ii)} \Longrightarrow \mathrm{(iii)}]$ is obvious.

$[ \mathrm{(iii)} \Longrightarrow \mathrm{(i)}]$ By Lemma \ref{lem81} it
suffices to check that the inclusion $`` \subset"$ holds. Let $x^* \in
\partial p (x)$. We have $x \in (Mp)(x^*)$. Since robust duality holds at $%
x^*$, Theorem \ref{thm31} (with $\varepsilon = 0$) says that $x \in \mathcal{%
A} (x^*) $. Thus, $x^* \in \mathcal{A} ^{-1} (x)$, and, by Lemma \ref{lem61}%
, $x^* \in C(x)$. \qed 

In the deterministic and the non-parametric cases, we get the next results
from Theorem \ref{thm81}.

\begin{corollary}
\label{corol81} Let $F:X\times Y\rightarrow \mathbb{R}_{\infty }$, $%
p=F(\cdot ,0_{Y})$, and $x\in p^{-1}(\mathbb{R})$. The next statements are
equivalent:

$\mathrm{(i)} $\ \ \ $\partial p (x) = \bigcap\limits_{\eta > 0 } \mathrm{%
proj}_{X^\ast} (\partial^\eta F)(x, 0_Y)$,

$\mathrm{(ii)} $\ \ \ $\min\limits_{z \in X} \Big\{ F(z, 0_Y) - \langle x^*,
x \rangle \Big\} = \sup\limits_{y^* \in Y^\ast} - F^\ast (x^*, y^*) , \
\forall x^* \in \partial p (x)$,

$\mathrm{(iii)} $ \ $\inf\limits_{z \in X} \Big\{ F(z, 0_Y) - \langle x^*, x
\rangle \big\} = \sup\limits_{y^* \in Y^\ast} - F^\ast (x^*, y^*) , \
\forall x^* \in \partial p (x)$.
\end{corollary}

\textit{Proof.} Let $F_{u}=F:X\times Y\rightarrow \mathbb{R}_{\infty }$ and $%
p=F(\cdot ,0_{Y})$. We then have 
\begin{equation*}
C(x)=\bigcap\limits_{\eta >0}\mathrm{proj}_{X^{\ast }}(\partial ^{\eta
}F)(x,0_{Y}),\ \forall x\in X,
\end{equation*}%
(see Corollary \ref{corol61}) and the conclusion follows directly from
Theorem \ref{thm81}. \qed

\begin{corollary}
\label{corol82} Let $(f_{u})_{u\in U}\subset \mathbb{R}_{\infty }^{X}$, $%
p=\sup\limits_{u\in U}f_{u}$, $x\in p^{-1}(\mathbb{R})$. The next statements
are equivalent:

$\mathrm{(i)} $\ \ $\partial \left( \sup\limits_{u \in U} f_u\right) (x) =
C(x)$,

$\mathrm{(ii)} $\ \ $\max\limits_{z \in X} \Big\{  \langle x^*, z \rangle -
p(z) \Big\} = \inf\limits_{u \in U} f_u^\ast (x^*), \ \forall x^* \in
\partial p (x)$,

$\mathrm{(iii)} $ \ $\left( \sup\limits_{u \in U} f_u\right) ^\ast (x^*) =
\inf\limits_{u \in U } f_u^\ast (x^*), \ \forall x^* \in \partial p (x)$,

where 
\begin{equation}
C(x)=\bigcap\limits_{\eta >0}\bigcup\limits_{\QATOP{\scriptstyle{\varepsilon
_{1}}+{\varepsilon _{2}}=\eta \hfill }{\scriptstyle{\varepsilon _{1}}%
\geqslant 0,{\varepsilon _{2}}\geqslant 0\hfill }}\bigcup\limits_{u\in
I^{\varepsilon _{1}}(x)}\mathrm{proj}^u_{X^{\ast }}(\partial ^{\varepsilon
_{2}}f_{u})(x),\forall x\in X.  \label{83}
\end{equation}
\end{corollary}

\textit{Proof.} Let $F_u(x, y_u) = f_u(x)$. Then it is easy to see that in
this case, $C(x) $ can be expressed as in \eqref{83}. The conclusion now
follows from Theorem \ref{thm81}. \qed

\medskip

Let us come back to the general case and consider the most simple
subdifferential formula one can expect for the robust objective function $%
p=\sup\limits_{u\in U}F_{u}(\cdot ,0_{u})$: 
\begin{equation}
\partial p(x)=\displaystyle\bigcup_{u\in I(x)}\mathrm{proj}^u_{X^{\ast
}}\left( \partial F_{u}\right) (x,0_{u}),  \label{84}
\end{equation}%
where the set of active indexes at $x$, $I(x)$, is defined by \eqref{I-x}.

In Case 3 \ we have 
\begin{equation*}
p(x)=\left\{ 
\begin{array}{ll}
f(x), & \mathrm{if}\quad H_{u}(x)\in -S_{u},\forall u\in U, \\ 
+\infty , & \mathrm{else},%
\end{array}%
\right.
\end{equation*}%
$I(x)=U$ for each $x\in p^{-1}(\mathbb{R})$, and \eqref{84} writes 
\begin{equation*}
\partial p(x)=\bigcup\limits_{\QATOP{u\in U,\ z_{u}^{\ast }\in S_{u}^{+}}{%
\langle z_{u}^{\ast },H_{u}(x)\rangle =0}}\partial (f+z_{u}^{\ast }\circ
H_{u})(x),
\end{equation*}%
which has been called \textit{Basic Robust Subdifferential Condition} (BRSC)
in \cite{BJL13} (see \cite[page 307]{HU-Lem1} for the deterministic case).
More generally, let us introduce the following terminology:

\begin{definition}
Given $F_{u}:X\times Y_{u}\rightarrow \mathbb{R}_{\infty }$ for each $u\in U$%
, and $p=\sup\limits_{u\in U}F_{u}(\cdot ,0_{u})$, we will say that \textbf{%
Basic Robust Subdifferential Condition holds at a point } $x\in p^{-1}(%
\mathbb{R})$ if \eqref{84} is satisfied, that is $\partial p(x)=D(x)$.
\end{definition}

Recall that, in Example \ref{Example1}, $p\left( x\right) =\left\langle
c^{\ast },x\right\rangle +\mathrm{i}_{A}\left( x\right) ,$ where $A=p^{-1}(%
\mathbb{R})$ is the feasible set of the linear system. So, given $x\in A,$ $%
\partial p(x)$ is the sum of $c^{\ast }$ with the normal cone of $A$ at $x,$
i.e., Basic Robust Subdifferential Condition (at $x$) asserts that such a
cone can be expressed in some prescribed way.

\begin{theorem}
\label{thm82} Let $x\in p^{-1}(\mathbb{R})$. The next statements are
equivalent:

$\mathrm{(i)}$\ \ Basic Robust Subdifferential Condition holds at $x$,

$\mathrm{(ii)}$\ Min-max robust duality holds at each $x^{\ast }\in \partial
p(x)$,

$\mathrm{(iii)}$ Strong robust duality holds at each $x^{\ast }\in \partial
p(x)$.
\end{theorem}

\noindent \textit{Proof.} $[\mathrm{(i)} \Longrightarrow \mathrm{(ii)}]$ Let 
$x^* \in \partial p (x)$. We have $x^* \in D(x)$ and by \eqref{81} there
exist $u \in I(x)$ (i.e., $p(x) = F_u(x, 0_u)$), $y_u^*\in Y_u^\ast$, such
that $(x^*, y_u^*) \in \partial F_u)(x, 0_u)$. Then, 
\begin{eqnarray*}
p^\ast (x^*) &\geq& \langle x^*, x\rangle - p(x) = \langle x^*, x \rangle -
F_u(x, 0_u) = F_u^\ast (x^*, y_u^*) \\
&\geq& q(x^*) \geq p^*(x^*) .
\end{eqnarray*}
It follows that 
\begin{equation*}
\max\limits_{z \in X} \{ \langle x^*, z\rangle - p(z) \} = \langle x^*,
x\rangle - p(x) = F_u^* (x^*, y_u^*) = q(x^*),
\end{equation*}
and min-max robust duality holds at $x^*$.

$[\mathrm{(ii)} \Longrightarrow \mathrm{(iii)}]$ is obvious.

$[\mathrm{(iii)} \Longrightarrow \mathrm{(i)}]$ By Lemma \ref{lem81}, it
suffices to check that $\partial p (x) \subset D(x)$. Let $x^* \in \partial
p (x)$. We have $x \in (M p) (x^*)$. Since strong robust duality holds at $%
x^*$, Theorem \ref{thm41} says that there exist $u \in U$, $y_u^* \in
Y_u^\ast$ such that $x \in B^0_{(u, y_u^*)} (x^*)$, that means (see %
\eqref{52}) 
\begin{equation*}
(x, 0_u) \in (MF_u)((x^*, y_u^*) , \ (x^*, y_u^*) \in \partial F_u)(x, 0_u),
\end{equation*}
and by \eqref{81}, $x^* \in D(x)$. \qed 

As usual, Theorem \ref{thm82} gives us corresponding results for the two
extreme cases: non-uncertainty and non-perturbation cases.

\begin{corollary}
\label{corol83} Let $F:X\times Y\rightarrow \mathbb{R}_{\infty }$, $%
p=F(\cdot ,0_{Y})$, and $x\in p^{-1}(\mathbb{R})$. The next statements are
equivalent:

$\mathrm{(i)} $\ \ $\partial p (x) = \mathrm{proj}_{X^\ast} (\partial F) (x,
0_Y)$,

$\mathrm{(ii)} $\ \ $\max\limits_{z \in X} \Big\{  \langle x^*, z \rangle -
F(z, 0_Y) \Big\} = \min\limits_{y^*\in Y^\ast} F^\ast (x^*, y^*), \ \forall
x^* \in \partial p (x)$,

$\mathrm{(iii)} $ \ $p ^\ast (x^*) = \min\limits_{y^* \in Y^\ast } F^\ast
(x^*, y^*), \ \forall x^* \in \partial p (x)$.
\end{corollary}

\textit{Proof.} In this case we have, by \eqref{72}, $D(x)=\mathrm{proj}%
_{X^{\ast }}(\partial F)(x,0_{Y}) $ and the conclusion is a direct
consequence of Theorem \ref{thm82}. \qed 

\begin{corollary}
\label{corol84} Let $(f_{u})_{u\in U}\subset \mathbb{R}_{\infty }^{X}$, $%
p=\sup\limits_{u\in U}f_{u}$, $x\in p^{-1}(\mathbb{R})$. The next statements
are equivalent:

$\mathrm{(i)} $ \ $\partial p (x) = \bigcup\limits_{u \in I(x)} \partial f_u
(x) $,

$\mathrm{(ii)} $\ \ $\max\limits_{z \in X} \Big\{  \langle x^*, z \rangle -
p(z) \Big\} = \min\limits_{u \in U} f_u^\ast (x^*), \ \forall x^* \in
\partial p (x)$,

$\mathrm{(iii)} $ \ $(\sup\limits_{u \in U} f_u) ^\ast (x^*) =
\min\limits_{y^* \in Y^\ast } f_u^\ast (x^*), \ \forall x^* \in \partial p
(x)$.
\end{corollary}

\textit{Proof.} In this non-parametric case, let $F_u(x, y_u) = f_u(x)$, $p
= \sup\limits_{u \in U} f_u$. We have 
\begin{equation*}
D(x) = \bigcup_{u \in I(x)} \partial f_u (x), \ I(x) = \{ u \in U\, :\,
f_u(x) = p(x) \in \mathbb{R}\}
\end{equation*}
and the Theorem \ref{thm82} applies. \qed 

\textbf{Acknowledgements} This research was supported by the National
Foundation for Science \& Technology Development (NAFOSTED) from Vietnam,
Project 101.01-2015.27 \textit{Generalizations of Farkas lemma with
applications to optimization}, by the Ministry of Economy and
Competitiveness of Spain and the European Regional Development Fund (ERDF)
of the European Commission, Project MTM2014-59179-C2-1-P, and by the
Australian Research Council, Project DP160100854. \textbf{\ }


\begin{thebibliography}{99}
\bibitem{BEN09} Ben-Tal, A., El Ghaoui, L., Nemirovski, A.: Robust
Optimization, Princeton U.P., Princeton (2009)

\bibitem{BS06} Bertsimas, D., Sim, M.: Tractable approximations to robust
conic optimization problems. Math. Programing \textbf{107B}, 5-36 (2006)

\bibitem{Bot-book1} Bot, R.I.: Conjugate duality in convex optimization.
Springer, Berlin (2010)

\bibitem{BJL13} Bot, R.I., Jeyakumar, V., Li, G.Y.: Robust duality in
parametric convex optimization. Set-Valued Var. Anal. \textbf{21}, 177-189
(2013)

\bibitem{BuJeWu06} Burachick, R.S., Jeyakumar, V., Wu, Z.-Y.: Necessary and
sufficient condition for stable conjugate duality. Nonlinear Anal.\textit{\ }%
\textbf{64},\textbf{\ }1998-2006 (2006)

\bibitem{Chu66} Chu, Y.C.: Generalization of some fundamental theorems on
linear inequalities. Acta Math Sinica \textbf{16}, 25-40 (1966)

\bibitem{DGLV17} Dinh, N., Goberna, M.A., L\'{o}pez, M.A., Volle, M.: A
unifying approach to robust convex infinite optimization duality. J. Optim.
Theory Appl. \textbf{174}, 650-685 (2017)

\bibitem{DMVV} Dinh, N., Mo, T.H., Vallet, G., Volle, M.: A unified approach
to robust Farkas-type results with applications to robust optimization
problems. SIAM J. Optim. \textbf{27}, 1075-1101 (2017)

\bibitem{DL17} Dinh, N., Long, D.H.: Complete characterizations of robust
strong duality for robust optimization problems. Vietnam J. Math., to
appear. DOI: 10.1007/s10013-018-0283-1

\bibitem{GL98} Goberna, M.A., L\'{o}pez, M.A.: Linear Semi-Infinite
Optimization. J. Wiley, Chichester (1998)

\bibitem{GL17} Goberna, M.A., L\'{o}pez, M.A.: Recent contributions to
linear semi-infinite optimization. 4OR-Q. J. Oper. Res. \textbf{15,} 221-264
(2017)

\bibitem{GLV14} Goberna, M.A., L\'{o}pez, M.A., Volle, M.: Primal attainment
in convex infinite optimization duality. J. Convex Anal. \textbf{21},
1043-1064 (2014)

\bibitem{Grad16} Grad, S.-M.: Closedness type regularity conditions in
convex optimization and beyond. Front. Appl. Math. Stat. \textbf{16},
September (2016) https://doi.org/10.3389/fams.2016.00014.

\bibitem{HU-Lem1} Hiriart-Urruty, J.-B., Lemarechal, C., Convex Analysis and
Minimization Algoritms I, Springer, Berlin (1993).

\bibitem{HMSV95} Hiriart-Urruty, J.-B., Moussaoui, M., Seeger, A., Volle,
M.: Subdifferential calculus without qualification conditions, using
approximate subdifferentials: a survey. Nonlinear Anal. \textbf{24},
1727-1754 (1995)

\bibitem{JL09} Jeyakumar, V. , Li, G.Y.: Stable zero duality gaps in convex
programming: complete dual characterisations with applications to
semidefinite programs. J. Math. Anal. Appl\textit{.} \textbf{360}, 156-167
(2009)

\bibitem{LJL11} Li, G.Y., Jeyakuma, V., Lee, G.M.: Robust conjugate duality
for convex optimization under uncertainty with application to data
classification. Nonlinear Anal. \textbf{74}, 2327-2341 (2011)

\bibitem{LR17} Lindsey, M., Rubinstein, Y.A.: Optimal transport via a
Monge-Amp\`{e}re optimization problem. SIAM J. Math. Anal. \textbf{49},
3073-3124 (2017)

\bibitem{LS07} L\'{o}pez, M.A., Still, G.: Semi-infinite programming.
European J. Oper. Res. \textbf{180}, 491-518 (2007)

\bibitem{Rock74} Rockafellar, R.T.: Conjugate Duality and Optimization. CBMS
Lecture Notes Series No. 162, SIAM (1974)

\bibitem{Vera17} Vera, J.R.: Geometric measures of convex sets and bounds on
problem sensitivity and robustness for conic linear optimization. Math.
Programming\textit{\ }\textbf{147A}, 47-79 (2014)

\bibitem{Volle95} Volle, M.: Caculus rules for global approximate minima and
applications to approximate subdifferential calculus. J. Global Optim. 
\textbf{5}, 131-157 (1994)

\bibitem{WSS15} Wang, Y., Shi, R., Shi, J.: Duality and robust duality for
special nonconvex homogeneous quadratic programming under certainty and
uncertainty environment. J. Global Optim. \textbf{62}, 643-659 (2015)

\bibitem{Za02} C. Z\u{a}linescu, Convex Analysis in General Vector Spaces.
World Scientific, River Edge, NJ (2002)
\end{thebibliography}
\end{document}